\documentclass[11pt,a4paper,reqno]{amsart}
\usepackage[T1]{fontenc}
\usepackage[utf8]{inputenc}
\usepackage{lmodern,amsmath,amssymb,amsthm,mathtools}
\usepackage[margin=26mm]{geometry}
\usepackage{microtype,cite,indentfirst,tikz,needspace}
\usepackage{placeins}
\usepackage[hidelinks]{hyperref}
\hypersetup{pdftitle={Polya inequalities on thin planar networks},
  pdfauthor={Kefeng Liu and Liangdi Zhang},
  pdfsubject={Spectral geometry of thin planar domains and Polya inequalities},
  pdfkeywords={spectral geometry, Polya conjecture, Laplace eigenvalues, thin domains}}
\allowdisplaybreaks[1]
\numberwithin{equation}{section}
\newtheorem{theorem}{Theorem}[section]
\newtheorem{proposition}[theorem]{Proposition}
\newtheorem{lemma}[theorem]{Lemma}
\newtheorem{corollary}[theorem]{Corollary}
\theoremstyle{remark}
\newtheorem{remark}[theorem]{Remark}
\newtheorem*{introremark}{Remark}
\newcommand{\R}{\mathbb R}
\newcommand{\cN}{c_{\mathrm N}}
\DeclareMathOperator{\dist}{dist}
\title[P\'olya inequalities on thin planar networks]{P\'olya inequalities on thin planar networks}
\author{Kefeng Liu}
\address{Mathematical Sciences Research Center, Chongqing University of
Technology, Chongqing 400054, China}
\email{kefengliu@cqut.edu.cn}
\author{Liangdi Zhang}
\address{Mathematical Sciences Research Center, Chongqing University of
Technology, Chongqing 400054, China}
\email{ldzhang91@163.com}
\date{November 13, 2026}
\makeatletter
\@ifundefined{subjclassname@2020}{%
  \@namedef{subjclassname@2020}{\textup{2020} Mathematics Subject Classification}}{}
\makeatother
\begin{document}
\begin{abstract}
We construct thin planar domains on which the Dirichlet
and the Neumann P\'olya inequalities both hold strictly. For the
neighborhoods of a given finite connected embedded polygonal graph
(branches and cycles are allowed) there is a single radius threshold
below which both inequalities hold for every eigenvalue index. We also
show that every nonempty bounded connected open set in the plane can be
approximated, in measure and in the Hausdorff distance, by smooth simply
connected domains satisfying both inequalities. Along the same sequence,
the first Dirichlet eigenvalue tends to infinity, while each Neumann
eigenvalue of fixed positive index tends to zero. In the special case of
rectangular frames our estimates lead to explicit sufficient conditions,
with no restriction on the aspect ratio.
\end{abstract}
\maketitle

\section{Introduction}\label{sec:introduction}

We study the Dirichlet and Neumann eigenvalues of thin neighborhoods
of finite planar graphs and their use in geometric approximation.
The networks may have branches and cycles, and the width condition
does not depend on the eigenvalue index.

Let $\Omega\subset\R^2$ be a bounded connected Lipschitz domain. Write
$\lambda_k(\Omega)$, $k\geq1$, and $\mu_k(\Omega)$, $k\geq0$, for its
Dirichlet and Neumann eigenvalues, listed in nondecreasing order and
counted with multiplicity, with $\mu_0(\Omega)=0$.
P\'olya's conjecture asserts that
\begin{equation}\label{eq:polya}
 |\Omega|\lambda_k(\Omega)\geq4\pi k,
 \qquad |\Omega|\mu_k(\Omega)\leq4\pi k
 \qquad(k\geq1).
\end{equation}
Here $|\Omega|$ denotes area. The constant $4\pi$ comes from
Weyl's law~\cite{hist_weyl1913}, so it cannot be improved uniformly in
$k$. The conjecture asks for an area bound on each eigenvalue, without
any convexity assumption or restriction on topology.
For spectral averages the corresponding sharp estimates are
known: we mention the Dirichlet bounds of Berezin~\cite{hist_berezin}
and of Li and Yau~\cite{hist_liyau}, the Neumann bounds of
Kr\"oger~\cite{hist_kroger}, and, more recently, the Riesz mean
asymptotics of Frank and Larson~\cite{hist_franklarsonriesz}. Such
averaged estimates, however, do not yield \eqref{eq:polya} index by
index. For an arbitrary bounded open set in the plane, the trace
identity of Speciel~\cite[Theorem~1]{speciel} gives
$\lambda_k\geq\frac23(4\pi k/|\Omega|)$ for every $k\geq1$, that is,
the conjectured bound with $2/3$ in place of the coefficient $1$.

P\'olya himself~\cite{polya1961} proved the Dirichlet inequality for
planar domains which tile the plane, and the Neumann inequality under
an additional regularity assumption on the tiling. This assumption was
later removed by Kellner~\cite{hist_kellner}.
Filonov~\cite{hist_filonovtiling} has given another proof for Neumann
tiling sets. (By a tiling we always mean that congruent copies of the
whole domain, with pairwise disjoint interiors, cover the plane up to
a set of area zero.) More recently, Filonov, Levitin, Polterovich and
Sher~\cite{hist_flpsballs} proved both
inequalities for disks and for circular sectors, together with the
Dirichlet inequality for balls in all dimensions. Refinements for balls
and cylinders can be found in the work of Guo, Miao, Wang and
Zhan~\cite{hist_guomiao}, and the Neumann problem on higher-dimensional
balls has been treated in two recent
preprints~\cite{hist_flpsneumann,hist_lineumann}.

Tang and Wang~\cite[Theorem~1.2]{tangwang} state the strict Dirichlet
inequality for every bounded simply connected planar domain.
They announce the corresponding strict Neumann inequality for bounded
simply connected Lipschitz domains in~\cite[Remark~1.3]{tangwang}.
The Neumann statement also appears in their public formalization
repository~\cite{tangwangNrepository}. These statements cover the
simply connected domains that occur below,
namely sufficiently thin neighborhoods of trees and the smooth
approximants of Theorem~\ref{thm:approximation}. We do not use either
result, and our arguments apply equally to fixed embedded networks
with cycles.

Thin domains are another natural source of examples. Already
Laptev~\cite{hist_laptev} observed that the Dirichlet P\'olya
inequality passes to a Cartesian product as soon as the factor which
satisfies it has dimension at least two. Freitas and
Salavessa~\cite[Lemma~1]{hist_freitassalavessa} showed how eventual
P\'olya bounds on a fixed domain carry over to sufficiently fine
partitions into congruent subdomains, and in this way obtained, among
other things, non-tiling thin cylinders. He and
Wang~\cite[Theorems~1.1 and~1.2]{hewang} proved both inequalities for
sufficiently thin Cartesian products. Jiang and
Lin~\cite[Theorems~1.9 and~1.15]{jianglin} obtained Dirichlet
inequalities for domains formed by removing suitable sets from a
strip-tiling domain. Their conditions concern both the removed sets
and the containing domain.

For networks with bends and branches, global separation of variables
is unavailable. We cut the domain into
rectangular strips and into fixed regions around the vertices, which we
call junctions. On the strips the mixed counting functions reduce to
sums over an interval, for which we use explicit estimates; on each
junction we use a Weyl remainder estimate for a fixed domain. Variational
bracketing then puts the pieces together in such a way that the exact
area term survives and the losses can be read off piece by piece.

The junctions may be general fixed regular domains
(see Section~\ref{sec:preliminaries}), and once the strips are thin
enough both strict inequalities follow. In the Neumann case the summed
counting bound leaves a range of low energies uncovered, and this range
is taken care of by low-energy trial functions. Trial functions which
vary only along the strips are also what allows the number of junctions
to grow in the approximation construction.

Let $\Gamma$ be a finite connected embedded polygonal graph with at
least one edge. In other words, its edges are closed line segments of positive length, and two
distinct edges can meet only at a common endpoint (in particular they
never overlap along a segment); polygonal graphs will always be
understood in this sense. We define the open $\varepsilon$-neighborhood by
\[
 \Omega_\varepsilon=\{x\in\R^2:\dist(x,\Gamma)<\varepsilon\}.
\]

\begin{theorem}\label{thm:graph}
For every finite connected embedded polygonal graph $\Gamma$ with at least
one edge, there exists
$\varepsilon_0=\varepsilon_0(\Gamma)>0$ such that
\begin{equation}\label{eq:graph-polya}
 |\Omega_\varepsilon|\lambda_k(\Omega_\varepsilon)>4\pi k,
 \qquad
 |\Omega_\varepsilon|\mu_k(\Omega_\varepsilon)<4\pi k
\end{equation}
for every $0<\varepsilon<\varepsilon_0$ and every $k\geq1$.
\end{theorem}

\begin{introremark}
The threshold $\varepsilon_0(\Gamma)$ depends on the graph but not on
the eigenvalue index, and nothing is assumed about the number of
branches or cycles. Since a sufficiently small regular neighborhood of
a graph with $M$ edges and $V$ vertices has $M-V+1$ holes, the theorem
covers in particular thin neighborhoods of polygonal cycles.

We note that these neighborhoods need not tile the plane. Take for
$\Gamma$ a single line segment, so that $\Omega_\varepsilon$ is a convex
stadium with semicircular ends. A tiling by congruent copies of it
would have to be locally finite, because the copies have a common
bounded diameter and each of them contains a disk of one and the same
positive radius, these disks being disjoint. Since the copies cover the
plane up to a set of area zero, the boundary of each copy would then
lie in the union of finitely many neighboring closed tiles. But two
convex tiles with disjoint interiors can meet only along a separating
line, and finitely many lines cannot cover a semicircular end; hence
the stadium does not tile the plane.
\end{introremark}

Let us compare this with the small-hole theorem of He, Jiang, Lin and
Wang~\cite[Theorem~1.1]{hjlw}, whose starting point is a fixed domain
satisfying an enhanced Neumann counting inequality. In
Theorem~\ref{thm:graph} the situation is different: the holes enclosed
by cycles may keep positive area while $|\Omega_\varepsilon|\to0$.
(We do borrow from \cite{hjlw} the interval estimate of
\cite[Lemma~2.6]{hjlw}, which is used in our Neumann argument.)

Thin neighborhoods of graphs have also been studied through convergence
to an operator on the graph. Freidlin and Wentzell~\cite{freidlinwentzell1993} considered
diffusion in narrow tubes. Kuchment and Zeng~\cite{kuchmentzeng2001}
and Rubinstein and Schatzman~\cite{rubinsteinschatzman2001} proved
convergence of the Neumann spectrum for domains shrinking onto a graph.
Exner and Post~\cite{exnerpost2005} dealt with graph-like manifolds,
and norm resolvent convergence was developed by
Post~\cite{post2006,post2012}. What these results describe is the
limiting graph operator. By themselves they do not control the sign
of the two-dimensional counting remainder at all energies (in
particular at energies where new transverse modes come in), and it is
exactly this sign that Theorem~\ref{thm:graph} needs.

Our second result is about geometric approximation. Since the area of
a neighborhood of one fixed graph tends to zero, we instead use
assemblies whose total strip length tends to infinity, so as to occupy
most of a given domain. That the spectrum can be unstable under perturbations of the
domain is classical (see e.g.\ Rauch and Taylor~\cite{rauchtaylor1975}),
and convergence theorems always need some extra control of the
domains or of the quadratic forms
\cite{mosco1994,kuwaeshioya2003,daners2003,burenkovlamberti2012}.
The construction below produces smooth simply connected approximants
which satisfy both strict P\'olya inequalities, although their Dirichlet
and Neumann spectra go in opposite directions. Here and below $d_H$
denotes the Hausdorff distance between nonempty compact sets.

\needspace{16\baselineskip}
\begin{theorem}[Approximation by smooth domains]\label{thm:approximation}
Let $U\subset\R^2$ be a nonempty bounded connected open set.
There is a sequence of simply connected domains $G_j$ with
$C^\infty$ boundary and $\overline{G_j}\subset U$ such that
\begin{equation}\label{eq:approximation-convergence}
 |U\setminus G_j|\longrightarrow0,
 \qquad d_H(\overline{G_j},\overline U)\longrightarrow0.
\end{equation}
For every $j\geq1$ and every $k\geq1$,
\begin{equation}\label{eq:approximation-polya}
 |G_j|\lambda_k(G_j)>4\pi k,
 \qquad |G_j|\mu_k(G_j)<4\pi k.
\end{equation}
The sequence can also be chosen so that
\begin{equation}\label{eq:approximation-spectrum}
 \lambda_1(G_j)\longrightarrow\infty,
 \qquad \mu_k(G_j)\longrightarrow0
 \quad\text{for every fixed }k\geq1.
\end{equation}
\end{theorem}

The approximation properties~\eqref{eq:approximation-convergence}, the two
inequalities in~\eqref{eq:approximation-polya} and the
limits~\eqref{eq:approximation-spectrum} all hold along one and the same
sequence. The sequence is constructed in Section~\ref{sec:approximation};
the number of junctions grows, and besides the assembly estimates we
need a space of longitudinal trial functions. The limits
in~\eqref{eq:approximation-spectrum} show that
convergence in measure and in the Hausdorff distance does not, by
itself, imply convergence of the spectra.

Rectangular frames provide an explicit application: the four corner
pieces have known spectra,
and so the assembly estimates yield concrete thickness conditions, for
either boundary condition and for every aspect ratio.

For comparison we recall what is known for annuli. Filonov, Levitin,
Polterovich and Sher~\cite{hist_flpsannuli} treated the Dirichlet
eigenvalues of circular annuli. Dai and Sun state strict inequalities
for all positive Neumann eigenvalues of annuli~\cite[Theorem~1.1]{daisunannuli},
and quantitative strict inequalities for annular
sectors~\cite[Theorems~1.1 and~1.2]{daisunsectors}. Li gives
strict Dirichlet inequalities for concentric spherical shells in every
dimension $d\geq2$~\cite[Theorem~1.1]{lishellD}, as well as the strict
Neumann inequalities~\cite[Theorem~1.1]{lishellN}. All of these
concern circular or spherical boundaries. For a rectangular frame we
keep the corners as they are and count the spectra of the rectangular
pieces. The small-hole theorem~\cite[Theorem~1.1]{hjlw} applies to
perforated rectangles too, but its limit (a small hole) is a different
one from the thin-frame limit studied here.

Theorem~\ref{thm:frames} allows
any rectangular aspect ratio and gives explicit sufficient thickness
conditions for the strict inequalities, under both boundary
conditions.

Let $a,b,h>0$. By a rectangular frame with inner side lengths $a,b$
and thickness $h$ we mean the open region lying between a rectangle
with side lengths $a+2h,b+2h$ and the concentric closed rectangle,
with parallel sides, of side lengths $a,b$.

\needspace{7\baselineskip}
\begin{theorem}\label{thm:frames}
Let $\Omega$ be a rectangular frame with inner side lengths $a,b>0$
and wall thickness $h>0$, with no restriction on $a/b$.
\begin{enumerate}
\item If $a+b\geq(18+6\sqrt2)h$, then
$|\Omega|\lambda_k(\Omega)>4\pi k$ for every $k\geq1$.
\item If $a+b\geq200h$, then
$|\Omega|\mu_k(\Omega)<4\pi k$ for every $k\geq1$.
\end{enumerate}
\end{theorem}

Note that both conditions involve only $(a+b)/h$; frames whose inner
side lengths are very different are therefore included.

The paper is organized as follows. In Section~\ref{sec:preliminaries}
we collect the variational comparisons we need and state our geometric
assumptions. Sections~\ref{sec:interval} and~\ref{sec:local} provide the
estimates for sums, rectangles and junctions, which are combined
into a counting theorem in Section~\ref{sec:assembly}.
Theorem~\ref{thm:graph} is proved in Section~\ref{sec:graph}, where we
also deal with smooth graph neighborhoods. Section~\ref{sec:chains}
contains the estimates for long chains which are needed for the proof
of Theorem~\ref{thm:approximation} in Section~\ref{sec:approximation}.
Theorem~\ref{thm:frames} is proved in Section~\ref{sec:frames}, and in
the final Section~\ref{sec:completion} we consider what happens when a
long rectangle is attached to a fixed core along an exposed straight
piece of its boundary.

\section{Preliminaries}\label{sec:preliminaries}

\subsection{Counting functions and variational comparisons}
Throughout, eigenvalues are counted with multiplicity. Whatever
boundary condition is under consideration, we denote the corresponding
eigenvalues by $\nu_j$ and put
\begin{equation}\label{eq:counts}
 N(E)=\#\{j:\nu_j<E\},\qquad
 \overline N(E)=\#\{j:\nu_j\leq E\}.
\end{equation}
The barred function also counts the eigenvalues equal to $E$.
A superscript $D$ or $N$ indicates a Dirichlet or a Neumann condition
on the whole boundary, and a subscript indicates the domain; note that
for $E>0$ the count $N_\Omega^N(E)$ always includes the constant
eigenfunction.

To prove the strict P\'olya inequalities, we establish the counting estimates
\begin{equation}\label{eq:counting-target}
 \overline N_\Omega^D(E)<\frac{|\Omega|E}{4\pi},\qquad
 N_\Omega^N(E)>\frac{|\Omega|E}{4\pi}
 \quad(E>0).
\end{equation}
Indeed, if we evaluate these at an eigenvalue, we get
\[
 k\leq\overline N_\Omega^D(\lambda_k)
 <\frac{|\Omega|\lambda_k}{4\pi},\qquad
 \frac{|\Omega|\mu_k}{4\pi}
 <N_\Omega^N(\mu_k)\leq k \quad(k\geq1).
\]
These inequalities also hold at multiple eigenvalues. Since the
spectra of rectangles look simpler in this form, we shall often write
$E=\pi^2t^2$ with $t>0$.

Recall that the quadratic form of the Laplacian is
\[
 q_\Omega(u)=\int_\Omega|\nabla u|^2,
 \qquad \|u\|_\Omega^2=\int_\Omega|u|^2.
\]
The form domain is $H_0^1(\Omega)$ in the Dirichlet case and
$H^1(\Omega)$ in the Neumann case; for mixed conditions one requires
the trace to vanish on the part of the boundary carrying the Dirichlet
condition.

Every domain whose spectrum we use is bounded and Lipschitz, so the
spectra are discrete. (The set $U$ to be approximated in
Theorem~\ref{thm:approximation} is not assumed to be Lipschitz, but we
never need its spectrum.) Comparisons will always be made through the
min--max principle for the Rayleigh quotients
$q_\Omega(u)/\|u\|_\Omega^2$; see for instance Safarov and
Vassiliev~\cite{hist_safarovbook}.

We use two consequences of the min--max principle. If we cut a domain
and allow the traces on the two sides of each
cut to be independent, the form domain becomes larger, so that the
eigenvalues go down and both counting functions go up. If instead we
impose zero traces on the cuts, the form domain becomes smaller: the
separated Dirichlet pieces embed into the original form domain by
extension by zero, the eigenvalues go up and the counting functions go
down.

For a finite partition, a function which is $H^1$ on each piece and
whose traces match across the interfaces belongs to $H^1$ of the
whole domain. Squared $L^2$ norms and energies add because the
interfaces have area zero.

A dilation by $\varepsilon>0$ multiplies the eigenvalues by
$\varepsilon^{-2}$, and so, for either count and for any given
boundary condition,
\[
 N_{\varepsilon J}(E)=N_J(\varepsilon^2E),\qquad
 \overline N_{\varepsilon J}(E)=\overline N_J(\varepsilon^2E).
\]
Finally, for nonempty compact sets $K_1,K_2\subset\R^2$ we use the
usual Hausdorff distance
\[
 d_H(K_1,K_2)=\max\left\{
 \sup_{x\in K_1}\dist(x,K_2),\sup_{y\in K_2}\dist(y,K_1)\right\}.
\]

\subsection{Junctions, ports, and admissible assemblies}\label{subsec:assemblies}
By a \emph{regular junction} we mean a bounded connected Lipschitz
domain whose boundary is made up of finitely many arcs, each smooth up
to its endpoints, which meet at interior angles lying in $(0,2\pi)$.
Typically the arcs will be line segments or circular arcs, and two
consecutive arcs may also join smoothly.

A \emph{port} is a relatively open straight segment contained in the
boundary of a junction. Different ports are disjoint, but they may
share an endpoint; since boundary conditions are understood in the
sense of traces, these finitely many endpoints play no role in the
form domain. We always assume that the part of the boundary lying
outside the closures of the ports contains an open arc of positive
length.

Fix junctions $J_v$, $1\leq v\leq V$, and pair their ports.
The two ports in each pair have the same length $b_e>0$.
They will be joined by a rectangular strip indexed by $1\leq e\leq M$.
We assume $V,M\geq1$.

At scale $\varepsilon$, the junctions are rigid images of
$\varepsilon J_v$. The strip $R_e$ has width $\varepsilon b_e$
and length $\ell_e(\varepsilon)>0$. Each end of width
$\varepsilon b_e$ coincides with the closure of its corresponding port.
A port need not fill a whole side of a polygon.

We call the assembly \emph{admissible} if the pieces have disjoint
interiors, any two of them meet only at paired ports or at endpoints
of ports, and
\[
 \Omega_\varepsilon=\operatorname{int}
 \left(\bigcup_v\overline{\varepsilon J_v}
       \ \cup\ \bigcup_e\overline{R_e}\right)
\]
is a connected Lipschitz domain (the rigid motions being understood
in this formula). We require the relative interior of every port to
lie inside $\Omega_\varepsilon$. Away from the interface endpoints,
the longitudinal strip sides and the nonport parts of the junction
boundaries must lie on $\partial\Omega_\varepsilon$.

Several pieces may share an interface endpoint. Whether that endpoint
lies on the physical boundary depends on the geometric realization.

The \emph{incidence graph} of an assembly has a vertex for each
junction and an edge for each strip, and the assembly is called a
\emph{chain} when this graph is a path. By the total strip length we
mean $\sum_e\ell_e$; junctions are not given any length. Recall also
that a regular neighborhood of an embedded graph is obtained by
replacing the vertices by topological disks and the edges by bands,
glued together according to the incidence relations and the planar
embedding of the graph.

Suppose that every junction closure is a disk with disjoint closed
ports, and that the pieces have no contacts other than the prescribed
attachments. If the incidence graph is a tree, then the closure of the
assembly is a disk.

To prove this, root the tree at one junction. Add the remaining
junctions one at a time, each together with the band joining it to its
predecessor. Attaching a band to a disk along one boundary interval
again gives a disk. Attaching the next disk to the other end of the
band has the same property. No other boundary points are identified,
so induction proves the claim.

Let $\xi_v$ denote the first eigenvalue on $J_v$ with Neumann
conditions on the ports and Dirichlet conditions on the rest of the
boundary. Then
\begin{equation}\label{eq:mixed-positive}
 \xi_v>0.
\end{equation}
For otherwise we could find $u_n$ in the mixed form domain with
\[
 \|u_n\|_{L^2(J_v)}=1,\qquad
 \|\nabla u_n\|_{L^2(J_v)}\longrightarrow0.
\]
By compactness a subsequence would converge, strongly in $L^2$
and weakly in $H^1$, to a constant. This constant has zero trace on the
Dirichlet arc and hence vanishes, which is absurd since its $L^2$ norm
is one.

\subsection{The fixed-domain Weyl estimate}
For a fixed regular junction $J$ we shall use the Weyl estimate with
remainder: there are constants $C_J\geq0$ and $E_J>0$ such that the
counting functions for pure Dirichlet and pure Neumann conditions
satisfy
\begin{equation}\label{eq:weyl-junction}
 \left|N_J^X(E)-\frac{|J|E}{4\pi}\right|\leq C_J\sqrt E,
 \qquad X\in\{D,N\},\quad E\geq E_J.
\end{equation}
This is the formulation for piecewise smooth domains given by He and
Wang~\cite[Section~2, equation~(2.6)]{hewang}, who refer to Seeley's
remainder estimates~\cite{seeley1978,seeley1980}; we only use it within
this piecewise smooth class, and the constants are allowed to depend
on $J$. Letting the energy decrease to $E$, one gets
\eqref{eq:weyl-junction} for $\overline N_J^X$ as well.

For suitable fixed domains, two-term asymptotics give the expected
signs at high energies. Ivrii~\cite{hist_ivrii} proved the two-term law
when the periodic billiard trajectories form a set of measure zero.
Melrose~\cite{hist_melrose} treated concave boundaries.
For the general theory, see Safarov and Vassiliev~\cite{hist_safarovbook}.
Riesz means and comparisons between
Dirichlet and Neumann eigenvalues give further information in the same
direction \cite{hist_franklarsonconvex,hist_filonovconvex,
hist_friedlander,hist_filonovcomparison,hist_safarovfilonov}, but none
of this is a substitute for a counting bound valid at every energy. In
our argument the remainder estimate \eqref{eq:weyl-junction} is used
only to control the junctions; the signs we need come from explicit
counts on rectangles.

Only remainders for pure Dirichlet and pure Neumann
conditions enter the counting argument. For the mixed junction problem,
we use only the positive first eigenvalue in \eqref{eq:mixed-positive}.
No mixed counting remainder is needed. The argument also requires
neither a signed second term nor constants uniform over varying
junction shapes.

\section{Two estimates for quarter-circle sums}\label{sec:interval}

On a rectangle, the counting estimates we need come down to sums of
the heights of a quarter circle at integer points. For $R\geq0$ let
\[
 f_R(x)=\sqrt{R^2-x^2}\quad(0\leq x\leq R),\qquad
 S(R)=\sum_{n=1}^{\lfloor R\rfloor}f_R(n).
\]
The quarter-disk area is $\int_0^R f_R(x)\,dx=\pi R^2/4$.
Define the deficit
\[
 D(R)=\frac{\pi R^2}{4}-S(R).
\]
(For $0\leq R<1$ the sum defining $S(R)$ is empty, so $S(R)=0$, and
when $R$ is a positive integer the last term $f_R(R)=\sqrt{R^2-R^2}$
vanishes.) We also put
\[
 \cN=1-\frac{\pi}{4}>0.
\]
The subscript $\mathrm N$ refers to Neumann conditions; $\cN$ will
keep this meaning throughout the paper. The next estimate is the case
$\gamma=1/2$ of He, Jiang, Lin and Wang~\cite[Lemma~2.6]{hjlw}, their
energy parameter being $R^2$, and we state it without proof.

\begin{lemma}[{\cite[Lemma~2.6]{hjlw}}]\label{lem:neumann-sum}
For every $R\geq 0$,
\begin{equation}\label{eq:neumann-sum}
 R+S(R)\geq \frac{\pi R^2}{4}+\cN R.
\end{equation}
Since equality holds for $R=1$, the constant $\cN$ is best possible.
\end{lemma}

A concavity argument for the Dirichlet sums can also be found in Jiang
and Lin~\cite[Lemma~2.1]{jianglin}. For our purposes it is important
that the deficit be bounded below by a linear function, with strict
inequality also at integer radii, and this is what the next lemma
provides.

\begin{lemma}\label{lem:dirichlet-sum}
For every $R\geq 1$,
\begin{equation}\label{eq:dirichlet-sum}
 D(R)>\frac{R}{3}.
\end{equation}
\end{lemma}

\begin{proof}
\smallskip
\noindent\emph{Step 1: $R\geq8/3$.}
Let $m=\lfloor R\rfloor$, $\theta=R-m$ and $f(x)=\sqrt{R^2-x^2}$. By
the strict concavity of $f$ we have
\[
 \int_{n-1}^{n}f(x)\,dx>
 \frac{f(n-1)+f(n)}2\qquad(1\leq n\leq m).
\]
On the remaining interval $[m,R]$ we have
\[
 \int_m^R f(x)\,dx\geq\frac{\theta f(m)}2.
\]
This also holds if $m=R$. Adding these inequalities over
$[0,1],\ldots,[m-1,m]$ and $[m,R]$, subtracting
$S(R)=\sum_{n=1}^m f(n)$ and recalling that $f(0)=R$, we find
\[
 D(R)>\frac{R-(1-\theta)\sqrt{R^2-m^2}}2
       \geq\frac R2-\sqrt{\frac{2R}{27}}.
\]
The second inequality uses
\[
 R^2-m^2\leq2R\theta,\qquad
 \max_{0\leq\theta\leq1}(1-\theta)\sqrt\theta
 =\frac{2}{3\sqrt3}.
\]
Since $\sqrt{2R/27}\leq R/6$ for $R\geq8/3$, this gives $D(R)>R/3$.

\smallskip
\noindent\emph{Step 2: $1\leq R\leq3$.}
Here we write $T=R^2$ and use $\pi>157/50$ together with the
elementary inequality
\[
 \sqrt{x}\leq\frac{x}{2q}+\frac q2
 \qquad(x\geq0,\ q>0),
\]
which follows from $(\sqrt{x}-q)^2\geq0$.

Suppose first that $1\leq T\leq4$, so that $S(\sqrt T)=\sqrt{T-1}$.
With $q=10/13$ and $q=5/4$ the last inequality gives
\[
 \sqrt{T-1}\leq\frac{13}{20}(T-1)+\frac5{13},
 \qquad \frac{\sqrt T}{3}\leq\frac{2T}{15}+\frac5{24}.
\]
Substituting these bounds into $D(\sqrt T)=\pi T/4-S(\sqrt T)$ gives
\[
 D(\sqrt T)-\frac{\sqrt T}{3}
 >\frac{157T}{200}-\frac{13}{20}(T-1)-\frac5{13}
       -\frac{2T}{15}-\frac5{24}
 =\frac{T}{600}+\frac{89}{1560}>0.
\]

If $4\leq T\leq9$, then
$S(\sqrt T)=\sqrt{T-1}+\sqrt{T-4}$, and the choices
$q=2,15/14,9/4$ give, respectively,
\[
 \sqrt{T-1}\leq\frac T4+\frac34,\qquad
 \sqrt{T-4}\leq\frac7{15}(T-4)+\frac{15}{28},\qquad
 \frac{\sqrt T}{3}\leq\frac{2T}{27}+\frac38.
\]
Hence
\begin{align*}
 D(\sqrt T)-\frac{\sqrt T}{3}
 &>\frac{157T}{200}-\frac T4-\frac34
    -\frac7{15}(T-4)-\frac{15}{28}-\frac{2T}{27}-\frac38\\
 &=-\frac{31T}{5400}+\frac{173}{840}
 \geq\frac{27}{175}>0.
\end{align*}
The formulas for $S(\sqrt T)$ remain valid at $T=4$ and $T=9$,
because the extra endpoint term is zero there. Thus $D(R)>R/3$ holds
both for $1\leq R\leq3$ and for $R\geq8/3$, that is, for all $R\geq1$,
and \eqref{eq:dirichlet-sum} is proved.
\end{proof}

\section{Rectangles and fixed junctions}\label{sec:local}

\subsection{Mixed conditions on a rectangle}
Let $R_{\ell,w}=(0,\ell)\times(0,w)$, with $\ell,w>0$.
The sides of length $\ell$ will be called longitudinal and those of
length $w$ the ends; we do not assume that $\ell\geq w$ or the
opposite. We need two counting functions:
\begin{itemize}
\item $\mathcal D_{\ell,w}(E)$ counts eigenvalues at most $E$ with
Dirichlet conditions on the longitudinal sides and Neumann conditions
at the ends;
\item $\mathcal N_{\ell,w}(E)$ counts eigenvalues below $E$ with
Neumann conditions on the longitudinal sides and Dirichlet conditions
at the ends.
\end{itemize}
The conditions at the ends are the ones that come from cutting the
domain at the ports.

\begin{lemma}\label{lem:rectangle}
For every $t>0$ (including the values for which $\pi^2t^2$ is an
eigenvalue of the rectangle) we have
\begin{align}
 \mathcal D_{\ell,w}(\pi^2t^2)&=0
 &&\text{if }wt<1,\label{eq:rectangle-gap}\\
 \frac{\pi\ell wt^2}{4}-\mathcal D_{\ell,w}(\pi^2t^2)
 &>\left(\frac{\ell}{3}-w\right)t
 &&\text{if }wt\geq 1,\label{eq:rectangle-dirichlet}\\
 \mathcal N_{\ell,w}(\pi^2t^2)-\frac{\pi\ell wt^2}{4}
 &\geq (\cN\ell-w)t-1.
 &&\label{eq:rectangle-neumann}
\end{align}
\end{lemma}

\begin{proof}
Throughout the proof, $R=wt$ and $z_n=(\ell/w)\sqrt{R^2-n^2}$ for
$0\leq n\leq R$.

\smallskip
\noindent\emph{Step 1: Dirichlet conditions on the longitudinal sides.}
By separation of variables the eigenvalues are
\[
 \pi^2\left(\frac{m^2}{\ell^2}+\frac{n^2}{w^2}\right),
 \qquad m\geq0,\quad n\geq1.
\]
The smallest eigenvalue is $\pi^2/w^2$, attained for $(m,n)=(0,1)$.
If $wt<1$ then $\pi^2t^2<\pi^2/w^2$, so nothing is counted; this is
\eqref{eq:rectangle-gap}.

For a given $n\leq R$, the inclusive count allows exactly
$m=0,\ldots,\lfloor z_n\rfloor$, and therefore
\[
 \mathcal D_{\ell,w}(\pi^2t^2)
 =\sum_{n=1}^{\lfloor R\rfloor}(\lfloor z_n\rfloor+1)
 \leq\frac\ell w S(R)+\lfloor R\rfloor.
\]
Here we only used $\lfloor z_n\rfloor\leq z_n$ in each row and kept
all the $\lfloor R\rfloor$ terms $+1$; this matters when $R$ is an
integer, since the row $n=R$ then still contains the mode $m=0$. If
$R\geq1$, it follows from Lemma~\ref{lem:dirichlet-sum} that
\[
 \frac{\pi\ell wt^2}{4}-\mathcal D_{\ell,w}(\pi^2t^2)
 \geq\frac\ell wD(R)-\lfloor R\rfloor
 >\frac{\ell R}{3w}-R
 =\left(\frac\ell3-w\right)t.
\]
This proves \eqref{eq:rectangle-dirichlet}.

\smallskip
\noindent\emph{Step 2: Neumann conditions on the longitudinal sides.}
The eigenvalues are given by the same formula, except that now
$m\geq1$ and $n\geq0$. For $z\geq0$ let
\[
 p(z)=\#\{m\in\mathbb Z:m\geq1,\ m<z\}
     =\max\{0,\lceil z\rceil-1\}\geq z-1.
\]
Counting strictly in each row, we obtain
\begin{align*}
 \mathcal N_{\ell,w}(\pi^2t^2)
 &=\sum_{n=0}^{\lfloor R\rfloor}p(z_n)\\
 &\geq\frac\ell w\bigl(R+S(R)\bigr)-(\lfloor R\rfloor+1)\\
 &\geq\frac{\pi\ell wt^2}{4}+(\cN\ell-w)t-1.
\end{align*}
The last line follows from Lemma~\ref{lem:neumann-sum} and
$\lfloor R\rfloor\leq R$. When $z_n=0$ the row is empty and the
estimate for that row is just $0\geq-1$; so \eqref{eq:rectangle-neumann}
is also valid at the integer endpoints.
\end{proof}

\subsection{Bounds valid at every energy}
From the remainder estimate \eqref{eq:weyl-junction} we now derive
bounds which hold at all energies, not only at large ones. The
constant term in the Neumann bound is there because of the zero
eigenvalue.

\begin{lemma}\label{lem:junction}
For every fixed regular junction $J$, there are finite constants
$C_J^+,C_J^-\geq 0$ such that, for all $s>0$,
\begin{align}
 \overline N_J^{N}(\pi^2s^2)
 &\leq \frac{\pi|J|s^2}{4}+C_J^+s+1,
 \label{eq:junction-neumann}\\
 N_J^{D}(\pi^2s^2)
 &\geq \frac{\pi|J|s^2}{4}-C_J^-s.
 \label{eq:junction-dirichlet}
\end{align}
If $J$ is a rectangle with side lengths $a,b>0$, one can take both
constants equal to $a+b$.
\end{lemma}

\begin{proof}

\smallskip
\noindent\emph{Step 1: general $J$.}
Put $a_J=\pi|J|/4$. By \eqref{eq:weyl-junction} we can find
$s_1>0$ and $K\geq0$ such that, for $s\geq s_1$,
\[
 \overline N_J^N(\pi^2s^2)\leq a_Js^2+Ks,\qquad
 N_J^D(\pi^2s^2)\geq a_Js^2-Ks.
\]
Choose $0<s_0<s_1$ such that
$\pi^2s_0^2<\mu_1(J)$, and let
\[
 M_J=\overline N_J^N(\pi^2s_1^2),\qquad
 C_J^+=\max\left\{K,\frac{M_J-1}{s_0}\right\},\qquad
 C_J^-=\max\{K,a_Js_1\}.
\]
The count $M_J$ is finite. As $J$ is connected, the Neumann eigenvalue
zero is simple, so that $\overline N_J^N(\pi^2s^2)=1$ for
$0<s\leq s_0$ and \eqref{eq:junction-neumann} holds trivially in this
range. For $s_0\leq s\leq s_1$ we use monotonicity:
\[
 \overline N_J^N(\pi^2s^2)\leq M_J
 \leq1+C_J^+s\leq1+C_J^+s+a_Js^2.
\]
As for the Dirichlet bound \eqref{eq:junction-dirichlet}, when
$0<s\leq s_1$ it follows from
\[
 a_Js^2-C_J^-s\leq0\leq N_J^D(\pi^2s^2).
\]
For $s\geq s_1$, both \eqref{eq:junction-neumann} and
\eqref{eq:junction-dirichlet} are consequences of
\eqref{eq:weyl-junction}, because $C_J^\pm\geq K$. Thus both
estimates hold for every $s>0$.

\smallskip
\noindent\emph{Step 2: the Neumann bound for a rectangle.}
Let $J$ be a rectangle with side lengths $a,b$. Separating variables,
we get
\begin{align*}
 \overline N_{a,b}^N(\pi^2s^2)
 &=\sum_{n=0}^{\lfloor bs\rfloor}
 \left(\left\lfloor\frac ab\sqrt{(bs)^2-n^2}\right\rfloor+1\right)\\
 &\leq as+1+\frac ab S(bs)+\lfloor bs\rfloor\\
 &\leq\frac{\pi abs^2}{4}+(a+b)s+1.
\end{align*}
The last step uses $S(bs)\leq\int_0^{bs}\sqrt{(bs)^2-x^2}\,dx$,
the integrand being decreasing.

\smallskip
\noindent\emph{Step 3: the Dirichlet bound for a rectangle.}
Put $A=as$, $B=bs$, and
\[
 g(x)=B\sqrt{1-\frac{x^2}{A^2}}\quad (0\leq x\leq A),
 \qquad g(x)=0\quad (x>A).
\]
Since $g$ is nonincreasing,
\[
 \sum_{m=1}^{\lfloor A\rfloor}g(m)
 \geq\int_1^\infty g(x)\,dx
 \geq\frac{\pi AB}{4}-B.
\]
If in each row we replace $g(m)$ by the number of positive integers
strictly less than $g(m)$, we lose at most one; as there are
$\lfloor A\rfloor\leq A$ rows, we conclude that
\[
 N_{a,b}^{D}(\pi^2s^2)
 \geq \frac{\pi abs^2}{4}-(a+b)s.\qedhere
\]
\end{proof}

\section{The assembly theorem}\label{sec:assembly}

In this section we consider an admissible assembly in the sense of
Section~\ref{subsec:assemblies}, with reference junctions $J_v$, strip
widths $\varepsilon b_e$ and strip lengths $\ell_e=\ell_e(\varepsilon)$.
Recall from \eqref{eq:mixed-positive} that the mixed eigenvalues
$\xi_v$ of the reference junctions are positive. For the $e$th strip we
abbreviate the mixed counts as
\[
 \mathcal D_e(E)=\mathcal D_{\ell_e,\varepsilon b_e}(E),\qquad
 \mathcal N_e(E)=\mathcal N_{\ell_e,\varepsilon b_e}(E).
\]

With constants $C_{J_v}^{\pm}$ chosen as in Lemma~\ref{lem:junction},
we put
\begin{align}
 \tau&=\min\left\{\min_{1\leq v\leq V}\frac{\sqrt{\xi_v}}{\pi},
                    \min_{1\leq e\leq M}\frac{1}{b_e}\right\}>0,
 \label{eq:assembly-tau}\\
 d_e&=\min\left\{\frac13,\frac{\pi b_e\tau}{4}\right\},
 \qquad
 A_\varepsilon=\varepsilon\sum_e b_e\ell_e
             +\varepsilon^2\sum_v|J_v|,
 \label{eq:assembly-area}\\
 D_\varepsilon&=\sum_e d_e\ell_e
       -\varepsilon\left(\sum_e b_e+\sum_vC_{J_v}^+
                                    +\frac{V}{\tau}\right),
 \label{eq:assembly-D}\\
 B_\varepsilon&=\cN\sum_e\ell_e
       -\varepsilon\left(\sum_e b_e+\sum_vC_{J_v}^-\right).
 \label{eq:assembly-B}
\end{align}
Here $A_\varepsilon$ is the exact area, and $\tau$ controls the mixed
spectral gap. The coefficients $D_\varepsilon$ and $B_\varepsilon$
enter the Dirichlet deficit \eqref{eq:assembly-dirichlet} and the
Neumann excess \eqref{eq:assembly-neumann}, respectively.
The subtracted terms in \eqref{eq:assembly-D} and
\eqref{eq:assembly-B} account for strip and junction losses;
the term $V/\tau$ absorbs the constant terms in the junction counts.

The constants $\tau$, $d_e$ and $C_{J_v}^{\pm}$ depend only on the
reference pieces and their ports, not on the strip lengths or on
$\varepsilon$.

\needspace{9\baselineskip}
The main result of this section turns the local estimates into
counting bounds for the assembled domain.

\begin{theorem}\label{thm:assembly}
Let $\Omega_\varepsilon$ be an admissible assembly.
\begin{enumerate}
\item The Dirichlet spectrum satisfies
\begin{equation}\label{eq:assembly-gap}
 \lambda_1(\Omega_\varepsilon)
 \geq\frac{\pi^2\tau^2}{\varepsilon^2},
\end{equation}
\begin{equation}\label{eq:assembly-dirichlet}
 \overline N_{\Omega_\varepsilon}^{D}(\pi^2t^2)
 <\frac{\pi A_\varepsilon t^2}{4}-D_\varepsilon t
 \qquad\left(t\geq\frac{\tau}{\varepsilon}\right).
\end{equation}
If $D_\varepsilon\geq 0$, then
$A_\varepsilon\lambda_k(\Omega_\varepsilon)>4\pi k$ for every $k\geq 1$.
\item The Neumann counting function satisfies
\begin{equation}\label{eq:assembly-neumann}
 N_{\Omega_\varepsilon}^{N}(\pi^2t^2)
 \geq\max\left\{1,\frac{\pi A_\varepsilon t^2}{4}
                       +B_\varepsilon t-M\right\}
 \qquad(t>0).
\end{equation}
If
\begin{equation}\label{eq:assembly-neumann-condition}
 B_\varepsilon>\frac{M}{2}\sqrt{\pi A_\varepsilon},
\end{equation}
then $A_\varepsilon\mu_k(\Omega_\varepsilon)<4\pi k$ for every $k\geq 1$.
\end{enumerate}
\end{theorem}

\begin{lemma}[Dirichlet localization]\label{lem:opt-dirichlet-localization}
For an admissible assembly,
\[
 \lambda_1(\Omega_\varepsilon)\geq\frac{\pi^2\tau^2}{\varepsilon^2}.
\]
Moreover, for $t\geq\tau/\varepsilon$,
\[
 \frac{\pi A_\varepsilon t^2}{4}
 -\overline N_{\Omega_\varepsilon}^{D}(\pi^2t^2)>D_\varepsilon t.
\]
\end{lemma}

\begin{proof}
\smallskip
\noindent\emph{The spectral gap.}
Take $u\in H_0^1(\Omega_\varepsilon)$ and restrict it to the junctions
and to the rectangles; the restrictions still have zero trace on
$\partial\Omega_\varepsilon$. If we drop the requirement that the
traces of the two adjacent pieces agree at each port, the form domain
gets larger, and we are left with Neumann conditions on both sides of
every interface.

The resulting problems are the mixed junction problems defining
$\xi_v$ and the rectangles with Dirichlet conditions on the
longitudinal sides and Neumann conditions at the ends. Since the
squared $L^2$ norms and the Dirichlet energies are additive under
restriction, the min--max principle gives
\[
 \lambda_1(\Omega_\varepsilon)
 \geq\min\left\{\min_v\frac{\xi_v}{\varepsilon^2},
                 \min_e\frac{\pi^2}{\varepsilon^2b_e^2}\right\}
 =\frac{\pi^2\tau^2}{\varepsilon^2}.
\]

\smallskip
\noindent\emph{The counting deficit.}
To bound the counting function from above we go one step further and
also remove the Dirichlet conditions on the boundaries of the
junctions. This enlarges the form domain once more, while the problems
on the rectangles stay the same, and so
\begin{equation}\label{eq:dirichlet-bracketing}
 \overline N_{\Omega_\varepsilon}^{D}(\pi^2t^2)
 \leq\sum_v\overline N_{J_v}^{N}(\pi^2\varepsilon^2t^2)
       +\sum_e\mathcal D_e(\pi^2t^2).
\end{equation}

Now fix $t\geq\tau/\varepsilon$ and set
\[
 \delta_e(t)=\frac{\pi\varepsilon b_e\ell_et^2}{4}
                -\mathcal D_e(\pi^2t^2).
\]
By Lemma~\ref{lem:rectangle} and the definition of $d_e$,
\begin{align*}
 \varepsilon b_et\geq1:
 &\quad \delta_e(t)>
       \left(\frac{\ell_e}{3}-\varepsilon b_e\right)t
       \geq(d_e\ell_e-\varepsilon b_e)t,\\
 \varepsilon b_et<1:
 &\quad \delta_e(t)=\frac{\pi\varepsilon b_e\ell_et^2}{4}
       \geq\frac{\pi b_e\tau\ell_e}{4}t
       \geq d_e\ell_et
       >(d_e\ell_e-\varepsilon b_e)t.
\end{align*}

On the junctions we use the Neumann bound \eqref{eq:junction-neumann}
with $s=\varepsilon t$; summing over $v$, we get
\[
 \sum_v\left\{\frac{\pi\varepsilon^2|J_v|t^2}{4}
       -\overline N_{J_v}^{N}(\pi^2\varepsilon^2t^2)\right\}
 \geq-\varepsilon t\sum_vC_{J_v}^+-V.
\]

We now add up the $M$ deficits $\delta_e(t)$ of the rectangles and
the $V$ deficits of the junctions and use
\eqref{eq:dirichlet-bracketing}. The area terms add up to
$\pi A_\varepsilon t^2/4$, and since $t\geq\tau/\varepsilon$ we have
$V\leq\varepsilon(V/\tau)t$; hence
\begin{align*}
 \frac{\pi A_\varepsilon t^2}{4}
     -\overline N_{\Omega_\varepsilon}^{D}(\pi^2t^2)
 &>\left[\sum_ed_e\ell_e
       -\varepsilon\sum_eb_e-\varepsilon\sum_vC_{J_v}^+\right]t-V\\
 &\geq D_\varepsilon t.
\end{align*}
The inequality is strict because there is at least one strip.
\end{proof}

\begin{lemma}[The Neumann lower count]\label{lem:opt-neumann-count}
For an admissible assembly and every $t>0$,
\[
 N_{\Omega_\varepsilon}^{N}(\pi^2t^2)
 \geq\max\left\{1,\frac{\pi A_\varepsilon t^2}{4}
                      +B_\varepsilon t-M\right\}.
\]
\end{lemma}

\begin{proof}
This time we impose Dirichlet conditions on the boundaries of all
junctions and on the ends of the strips, and keep the Neumann
conditions on the longitudinal sides of the strips. A function in one
of these local form domains has zero trace on the interfaces, so its
extension by zero lies in $H^1(\Omega_\varepsilon)$, with the same
$L^2$ norm and the same Dirichlet energy.

For sums of such extensions the squared norms and the energies add
over the pieces, and the min--max principle yields
\begin{equation}\label{eq:neumann-bracketing}
 N_{\Omega_\varepsilon}^{N}(\pi^2t^2)
 \geq\sum_vN_{J_v}^{D}(\pi^2\varepsilon^2t^2)
       +\sum_e\mathcal N_e(\pi^2t^2).
\end{equation}
Inserting \eqref{eq:junction-dirichlet} and
\eqref{eq:rectangle-neumann}, we find
\begin{align*}
 N_{\Omega_\varepsilon}^{N}(\pi^2t^2)
 &\geq\frac{\pi A_\varepsilon t^2}{4}
      -\varepsilon t\sum_vC_{J_v}^-
      +\sum_e\bigl[(\cN\ell_e-\varepsilon b_e)t-1\bigr]\\
 &=\frac{\pi A_\varepsilon t^2}{4}+B_\varepsilon t-M.
\end{align*}
On the other hand, $u\equiv1$ is a Neumann eigenfunction on
$\Omega_\varepsilon$ with eigenvalue zero, so for $t>0$ we also have
$N_{\Omega_\varepsilon}^{N}(\pi^2t^2)\geq1$. The lemma follows by
taking the larger of the two bounds. (We do not claim that the
constant function can be added to the separated trial space; the two
bounds are used separately.)
\end{proof}

\begin{proof}[Proof of Theorem~\ref{thm:assembly}]
The estimates \eqref{eq:assembly-gap}--\eqref{eq:assembly-neumann}
were proved in Lemmas~\ref{lem:opt-dirichlet-localization}
and~\ref{lem:opt-neumann-count}, so it remains to deal with low
energies and then to evaluate the counting functions at eigenvalues.
Let
\[
 t_D=\frac\tau\varepsilon,\qquad
 t_N=\frac2{\sqrt{\pi A_\varepsilon}}.
\]
If $D_\varepsilon\geq0$, then
\[
 \begin{array}{ll}
 0<t<t_D:
  &\overline N_{\Omega_\varepsilon}^{D}(\pi^2t^2)=0
       <\pi A_\varepsilon t^2/4,\\[2pt]
 t\geq t_D:
  &\overline N_{\Omega_\varepsilon}^{D}(\pi^2t^2)
       <\pi A_\varepsilon t^2/4-D_\varepsilon t
       \leq\pi A_\varepsilon t^2/4.
 \end{array}
\]
If \eqref{eq:assembly-neumann-condition} holds, then
\[
 \begin{array}{ll}
 0<t<t_N:
  &N_{\Omega_\varepsilon}^{N}(\pi^2t^2)\geq1>
           \pi A_\varepsilon t^2/4,\\[2pt]
 t\geq t_N:
  &B_\varepsilon t-M\geq B_\varepsilon t_N-M>0.
 \end{array}
\]
For $t\geq t_N$, the strict Neumann counting bound follows from
the positivity of $B_\varepsilon t-M$ and \eqref{eq:assembly-neumann}.
So, under the respective hypotheses, the two bounds in
\eqref{eq:counting-target} hold for all $E>0$. Evaluating them at
eigenvalues (counted with multiplicity) we obtain
\[
 k\leq\overline N_{\Omega_\varepsilon}^{D}(\lambda_k)
       <\frac{A_\varepsilon\lambda_k}{4\pi},\qquad
 \frac{A_\varepsilon\mu_k}{4\pi}
       <N_{\Omega_\varepsilon}^{N}(\mu_k)\leq k.
\]
Multiplication by $4\pi$ gives the inequalities
$A_\varepsilon\lambda_k(\Omega_\varepsilon)>4\pi k$ and
$A_\varepsilon\mu_k(\Omega_\varepsilon)<4\pi k$, $k\geq1$, asserted in
the theorem.
\end{proof}

\begin{corollary}\label{cor:fixed-assembly}
Let the reference junctions, their ports, the incidence pattern and
the widths $b_e$ be fixed, and suppose that
\[
 \ell_e(\varepsilon)\longrightarrow L_e>0
 \qquad(\varepsilon\downarrow 0)
\]
for every strip. Assume that for all sufficiently small $\varepsilon$
there is an admissible geometric realization. Then there exists
$\varepsilon_0>0$ such that
\[
 |\Omega_\varepsilon|\lambda_k(\Omega_\varepsilon)>4\pi k,
 \qquad
 |\Omega_\varepsilon|\mu_k(\Omega_\varepsilon)<4\pi k.
\]
These inequalities hold for every $0<\varepsilon<\varepsilon_0$ and
every $k\geq 1$; a single $\varepsilon_0$ works for all eigenvalue indices.
\end{corollary}

\begin{proof}
All the reference constants in
\eqref{eq:assembly-tau}--\eqref{eq:assembly-B} are now fixed, and
\[
 D_\varepsilon\longrightarrow\sum_e d_eL_e>0,
 \qquad
 B_\varepsilon\longrightarrow\cN\sum_eL_e>0,
 \qquad
 A_\varepsilon\longrightarrow 0.
\]
Hence both sufficient conditions of Theorem~\ref{thm:assembly} are
satisfied once $\varepsilon$ is small enough. Since these conditions do
not involve the eigenvalue index, the uniformity in $k$ is automatic.
\end{proof}

We remark that the junctions need not be convex, nor rectangular, nor
simply connected; their shapes (holes included) are simply fixed before
we scale. The admissible thickness may depend on all these
data, and we make no claim of uniformity over a varying family of
junctions.

The enhanced Neumann counting inequality was introduced by He, Jiang,
Lin and Wang~\cite{hjlw}. As the next corollary shows, an inequality of
this form follows from \eqref{eq:assembly-neumann}, with constants that
can be written down explicitly in terms of the assembly.

\begin{corollary}[An enhanced Neumann bound]\label{cor:enhanced-neumann}
For an admissible assembly, write $A=|\Omega_\varepsilon|$ and
$B=B_\varepsilon$, and let $M\geq1$ be the number of strips.
Assume \eqref{eq:assembly-neumann-condition}. Set
\[
 \beta=\frac{\pi A}{4},\qquad
 r=\frac{\sqrt{B^2+4\beta(M+1)}-B}{2\beta},\qquad
 \eta=Br-M=1-\beta r^2.
\]
Then $M/B<r<1/\sqrt\beta$ and $\eta>0$. For every $E>0$,
\begin{equation}\label{eq:enhanced-neumann}
 N_{\Omega_\varepsilon}^N(E)\geq\frac{AE}{4\pi}
       +\max\left\{\eta,\frac{\eta}{\pi r}\sqrt E\right\}.
\end{equation}
\end{corollary}

\begin{proof}
Our assumption says that $B>M\sqrt\beta>0$. The polynomial
$F(t)=\beta t^2+Bt-(M+1)$ is strictly increasing for $t\geq0$, its only
positive zero is $r$, and
\[
 F(M/B)=\frac{\beta M^2}{B^2}-1<0,\qquad
 F(1/\sqrt\beta)=\frac{B}{\sqrt\beta}-M>0.
\]
These signs give the stated bounds on $r$. Moreover, $F(r)=0$ gives
$\eta=Br-M=1-\beta r^2$, and this is positive.

Next, by \eqref{eq:assembly-neumann} we have, for every $t>0$,
\[
 \Delta(t):=N_{\Omega_\varepsilon}^N(\pi^2t^2)-\beta t^2
 \geq\max\{1-\beta t^2,Bt-M\}.
\]
If $0<t\leq r$, then
\[
 \Delta(t)\geq1-\beta t^2\geq\eta
       =\max\{\eta,\eta t/r\}.
\]
If $t\geq r$, then
\[
 \Delta(t)\geq Bt-M
       =\frac{\eta t}{r}+M\left(\frac tr-1\right)
       \geq\frac{\eta t}{r}
       =\max\{\eta,\eta t/r\}.
\]
Putting $t=\sqrt E/\pi$ gives \eqref{eq:enhanced-neumann}.
The argument uses the strict count throughout and includes energies
equal to eigenvalues.
\end{proof}

Note that the corollary needs \eqref{eq:assembly-neumann-condition}.
In the approximation of Section~\ref{sec:approximation} we shall use
the subchain criterion of Proposition~\ref{prop:chain-criterion}
instead, and that criterion alone does not give the constants in
\eqref{eq:enhanced-neumann}.

\section{Polygonal graph neighborhoods}
\label{sec:graph}

We first decompose the distance neighborhood of a polygonal graph
exactly into strips and fixed reference junctions.

\begin{lemma}[Fixed reference junctions]\label{lem:opt-graph-decomposition}
Let $\Gamma$ be a fixed finite connected embedded polygonal graph with at
least one edge. For each vertex $v$, there are a fixed regular junction
$J_v$ and a number $R_v>1$ such that, for all sufficiently small
$\varepsilon>0$, its distance neighborhood has the admissible decomposition
\[
 \Omega_\varepsilon=
 \operatorname{int}\left(
   \bigcup_v(v+\varepsilon\overline{J_v})
       \cup\bigcup_e\overline{R_{e,\varepsilon}}\right).
\]
Every port has reference width $2$. For an edge $e=[v,w]$ of length $L_e$,
\[
 \operatorname{width}(R_{e,\varepsilon})=2\varepsilon,
 \qquad \ell_e(\varepsilon)=L_e-\varepsilon(R_v+R_w)>0.
\]
Each $J_v$ is bounded by finitely many line segments and circular
arcs, and has an open boundary arc which is not part of a port.
\end{lemma}

\begin{proof}
\smallskip
\noindent\emph{Step 1: fix the reference junctions.}
For a vertex $v$ and an edge $e$ incident to it, let $\omega_{v,e}$
be the unit vector pointing from $v$ along $e$. Define the ray star
and its unit neighborhood by
\[
 \Sigma_v=\bigcup_{e\ni v}\{r\omega_{v,e}:r\geq0\},\qquad
 U_v=\{x:\operatorname{dist}(x,\Sigma_v)<1\}.
\]
For each incident edge $e$, define the neighborhood of a single ray by
\[
 U_{v,e}=\{x\in\R^2:\operatorname{dist}
       (x,\{r\omega_{v,e}:r\geq0\})<1\}.
\]
If $e$ and $f$ are distinct edges at $v$, the intersection
$U_{v,e}\cap U_{v,f}$ is bounded. This also holds for opposite rays.
There are finitely many pairs, so choose $R_v>1$ with all these
intersections in $B(0,R_v-1)$. At a vertex with only one edge,
any $R_v>1$ will do.

We cut each outgoing strip at distance $R_v$ along the ray and let
$J_v$ be the bounded central piece. Every cut is a port of width $2$.
The numbers $R_v$ and the domains $J_v$ are now fixed.

The nonport boundary consists of line segments and circular arcs.
They meet smoothly or at angles in $(0,2\pi)$, and include an open
arc of positive length. Thus $J_v$ is regular and
$\xi_v>0$ by \eqref{eq:mixed-positive}.

The ray-star neighborhood is star-shaped with respect to the origin.
Truncating its disjoint outgoing strips preserves this property.
Hence $\overline{J_v}$ is a disk.

\smallskip
\noindent\emph{Step 2: place the pieces along the graph.}
Nonincident edges of the finite graph have positive separation.
Each vertex also has positive distance from every edge not containing it.
For small $\varepsilon$, nonincident pieces therefore remain disjoint.
Overlaps near a vertex $v$ lie inside
$v+\varepsilon\overline{J_v}$ by Step~1.

On an edge $e=[v,w]$, the junctions remove lengths
$\varepsilon R_v$ and $\varepsilon R_w$ from the two ends.
The remaining rectangle has
\[
 \operatorname{width}(R_{e,\varepsilon})=2\varepsilon,
 \qquad \ell_e(\varepsilon)
 =L_e-\varepsilon(R_v+R_w)>0.
\]
Its ends are precisely the two paired ports, and its longitudinal
sides lie on $\partial\Omega_\varepsilon$.

\smallskip
\noindent\emph{Step 3: verify admissibility.}
The junctions and rectangles cover $\Omega_\varepsilon$ and have
disjoint interiors. They meet only at paired ports or their endpoints.

The boundary consists of finitely many segments and circular arcs
meeting at admissible angles, so the union is Lipschitz.
It is connected because $\Gamma$ is connected.
These are the requirements of Section~\ref{subsec:assemblies}.
\end{proof}

\begin{proof}[Proof of Theorem~\ref{thm:graph}]
By Lemma~\ref{lem:opt-graph-decomposition}, $\Omega_\varepsilon$ is
an admissible assembly with
\[
 J_v\text{ and its ports fixed},\qquad b_e=2,\qquad
 \ell_e(\varepsilon)\longrightarrow L_e>0.
\]
Corollary~\ref{cor:fixed-assembly} gives a single
$\varepsilon_0(\Gamma)>0$ such that
\[
 |\Omega_\varepsilon|\lambda_k(\Omega_\varepsilon)>4\pi k,
 \qquad
 |\Omega_\varepsilon|\mu_k(\Omega_\varepsilon)<4\pi k
 \quad(0<\varepsilon<\varepsilon_0,\ k\geq1).
\]
\end{proof}

\needspace{11\baselineskip}
\begin{corollary}[Smooth graph neighborhoods]\label{cor:smooth-networks}
Let $\Gamma$ be a fixed finite connected embedded polygonal graph
with at least one edge and vertex set $\mathcal V$. There are constants $C,\varepsilon_0>0$ such that,
for every $0<\varepsilon<\varepsilon_0$, there is a bounded domain
$\widetilde\Omega_\varepsilon$ with $C^\infty$ boundary such that
\[
 \widetilde\Omega_\varepsilon=\Omega_\varepsilon
 \quad\text{outside }\bigcup_{v\in\mathcal V}B(v,C\varepsilon).
\]
Each $\widetilde\Omega_\varepsilon$ is a regular neighborhood of
$\Gamma$ and satisfies the strict Dirichlet and Neumann P\'olya
inequalities at every index. If $\Gamma$ has $V$ vertices and $M$ edges, then
$\widetilde\Omega_\varepsilon$ has exactly $M-V+1$ holes.
\end{corollary}

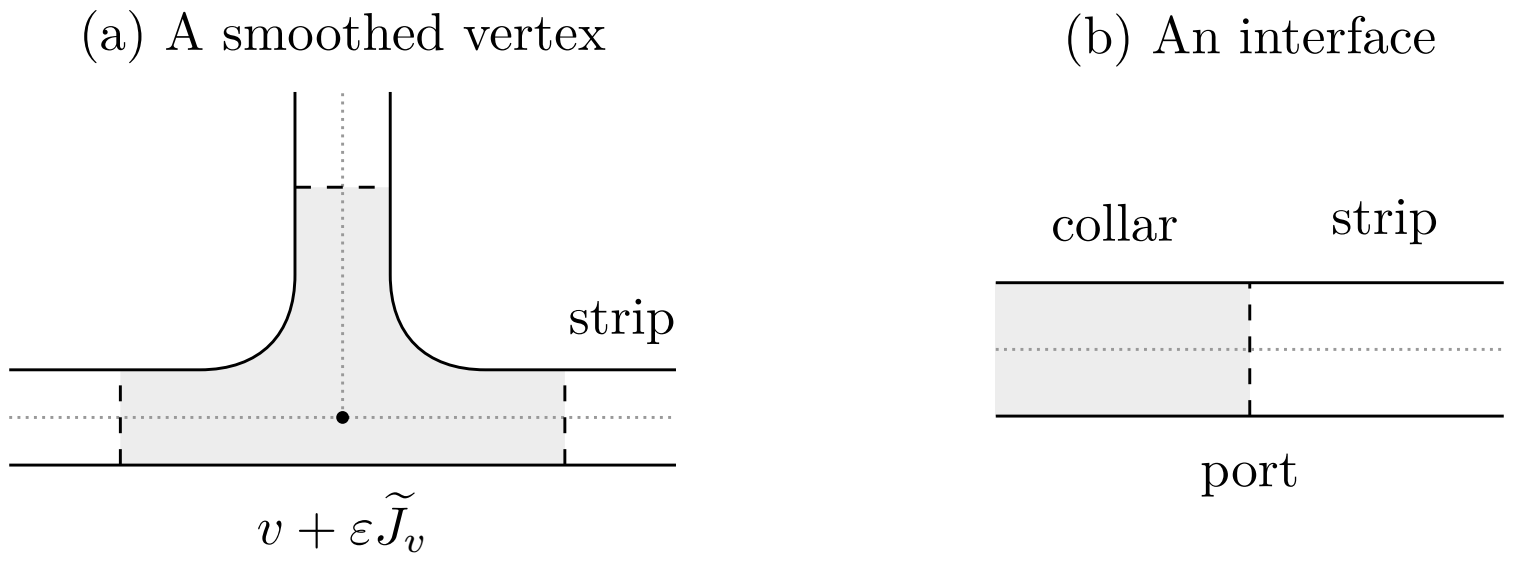
\begin{figure}[htbp]
\centering
\begin{tikzpicture}[x=1.05cm,y=1.05cm,line width=.55pt,font=\small]
 \fill[black!7]
 (-1.4,-.3)--(1.4,-.3)--(1.4,.3)--(.9,.3)
 .. controls (.52,.3) and (.3,.52) .. (.3,.9)
 --(.3,1.45)--(-.3,1.45)--(-.3,.9)
 .. controls (-.3,.52) and (-.52,.3) .. (-.9,.3)
 --(-1.4,.3)--cycle;
 \draw[black!40,densely dotted] (-2.1,0)--(2.1,0) (0,0)--(0,2.05);
 \draw (-2.1,-.3)--(2.1,-.3);
 \draw (2.1,.3)--(.9,.3)
 .. controls (.52,.3) and (.3,.52) .. (.3,.9)--(.3,2.05);
 \draw (-2.1,.3)--(-.9,.3)
 .. controls (-.52,.3) and (-.3,.52) .. (-.3,.9)--(-.3,2.05);
 \draw[dashed] (-1.4,-.3)--(-1.4,.3) (1.4,-.3)--(1.4,.3)
 (-.3,1.45)--(.3,1.45);
 \fill (0,0) circle (1.2pt);
 \node[below] at (0,-.34) {$v+\varepsilon\widetilde J_v$};
 \node[above] at (0,2.1) {(a) A smoothed vertex};
 \node at (1.76,.6) {strip};
 \begin{scope}[xshift=6cm,yshift=.45cm]
  \fill[black!7] (-1.6,-.42) rectangle (0,.42);
  \draw (-1.6,-.42)--(1.6,-.42) (-1.6,.42)--(1.6,.42);
  \draw[dashed] (0,-.42)--(0,.42);
  \draw[black!40,densely dotted] (-1.6,0)--(1.6,0);
  \node at (-.85,.8) {collar};
  \node at (.85,.8) {strip};
  \node[below] at (0,-.52) {port};
  \node[above] at (0,1.65) {(b) An interface};
 \end{scope}
\end{tikzpicture}
\caption{A schematic picture. The shaded junction is a scaled copy of
a fixed reference piece, and the dashed lines are the artificial
interfaces. Near each interface the boundary is straight, so that it
stays smooth after the strip is attached.}
\label{fig:smooth-junction}
\end{figure}

\begin{proof}
\smallskip
\noindent\emph{Step 1: smooth the reference junctions.}
Start from Lemma~\ref{lem:opt-graph-decomposition} and move each cut
a little further along its ray if needed. Each junction then contains
a rectangular neighborhood of every port, called a
\emph{rectangular collar}.

The collar has two straight boundary sides perpendicular to the port.
These sides continue along the attached strip.

Smooth the joins of the nonport boundary arcs in disjoint small
neighborhoods away from the collars. Choose the modifications so that
the closure remains a disk and the junction still contains the truncated
ray star. The collars remain unchanged.

This gives fixed domains $\widetilde J_v$ with smooth boundaries
away from the closures of their ports.

\smallskip
\noindent\emph{Step 2: assemble the smooth domain.}
Scale the new junctions by $\varepsilon$ and attach the rectangles.
The boundary is straight on both sides of every interface.
It is therefore $C^\infty$ across the interfaces
(Figure~\ref{fig:smooth-junction}).

The nonincident parts of the fixed graph have positive separation.
Thus no new contacts appear for small $\varepsilon$.
Since the reference modifications are fixed, they lie inside
$B(v,C\varepsilon)$ after scaling, for a single constant $C$.

\smallskip
\noindent\emph{Step 3: count the holes.}
The disk junctions and strips form a regular neighborhood of $\Gamma$.
Gluing $V$ disks and $M$ strips along $2M$ intervals gives
\[
 \chi(\widetilde\Omega_\varepsilon)=V+M-2M=V-M.
\]
Since $\chi=1-h$ for a connected planar domain with $h$ holes,
we get $h=M-V+1$.

\smallskip
\noindent\emph{Step 4: apply the assembly estimates.}
Every $\widetilde J_v$ has a nonport boundary arc of positive length,
so its mixed spectral gap is positive. Lemma~\ref{lem:junction}
gives finite counting constants for the fixed list of junctions.

The strip lengths satisfy $\ell_e(\varepsilon)\to L_e>0$.
Corollary~\ref{cor:fixed-assembly} therefore gives
\[
 |\widetilde\Omega_\varepsilon|
   \lambda_k(\widetilde\Omega_\varepsilon)>4\pi k,\qquad
 |\widetilde\Omega_\varepsilon|
   \mu_k(\widetilde\Omega_\varepsilon)<4\pi k
 \quad(k\geq1)
\]
for all sufficiently small $\varepsilon$.
The width threshold may depend on the smoothed junctions, which were
fixed before scaling.
\end{proof}

\begin{remark}
For instance, a polygonal cycle with two branches attached to it has
a smooth regular neighborhood with exactly one hole.
The admissible thickness may depend on the fixed graph's angles,
edge lengths, and separation of nonincident parts.
For Corollary~\ref{cor:smooth-networks}, it may also depend on the
smoothed junctions. The proof uses fixed templates and requires no
spectral continuity under smoothing.

A smooth curved edge cannot be decomposed exactly into rectangles.
Under polygonal approximation, the junctions change and the interface
losses accumulate. The next section treats certain growing chains;
it does not give a fixed-width result for arbitrary smooth tubes.
\end{remark}

\section{Long chains}\label{sec:chains}

The constant eigenfunction controls the Neumann count at very low
energies. On long chains, longitudinal trial functions extend this
control to a larger energy range. We need these bounds for the
approximation in Section~\ref{sec:approximation}.

In this section all strips have the same width $w=\varepsilon b$.
We keep the assumptions of Section~\ref{sec:assembly} and write
\[
 \Omega=\Omega_\varepsilon,\qquad A=|\Omega|,\qquad
 L=\sum_{e=1}^M\ell_e,\qquad B=B_\varepsilon,\qquad D=D_\varepsilon.
\]
Only the straight strips contribute to $L$; the bends lie in the
junctions. Chains are as in Section~\ref{subsec:assemblies}, and the
junctions may have holes.

\begin{lemma}[Longitudinal trial spaces]\label{lem:longitudinal}
If the assembly is a chain, then
\begin{align}
 \mu_k(\Omega)&\leq\frac{\pi^2k^2}{L^2}
 &&(k\geq0),\label{eq:chain-eigenvalues}\\
 N_\Omega^{N}(\pi^2t^2)&\geq\lceil Lt\rceil\geq Lt
 &&(t>0).\label{eq:chain-count}
\end{align}
More generally, suppose that the incidence graph contains a path
whose strips have total length $L_*>0$, and that none of the
intermediate junctions of this path meets a strip outside the path
(the two terminal junctions may have other attachments). Then
\begin{align}
 \mu_k(\Omega)&\leq\frac{\pi^2(k+1)^2}{L_*^2}
 &&(k\geq0),\label{eq:subchain-eigenvalues}\\
 N_\Omega^{N}(\pi^2t^2)&\geq
 \max\{1,\lceil L_*t\rceil-1\}
 \geq\max\{1,L_*t-1\}
 &&(t>0).\label{eq:subchain-count}
\end{align}
\end{lemma}

\begin{proof}
\smallskip
\noindent\emph{Step 1: the whole assembly is a chain.}
Order the rectangles along the path. Joining their longitudinal
coordinates gives $[0,L]$, with each junction represented by a point
$s_v$.

For $u\in H^1(0,L)$, define $Tu$ to equal $u(s)$ on each rectangle
cross section at $s$. Set $Tu=u(s_v)$ on the $v$th junction.

The traces match at every port, so $Tu\in H^1(\Omega)$.
The map $T$ is linear and injective. Its energy and mass satisfy
\begin{align}
 \int_\Omega|\nabla Tu|^2
 &=w\int_0^L|u'|^2,\label{eq:longitudinal-energy}\\
 \int_\Omega|Tu|^2
 &=w\int_0^L|u|^2+
   \varepsilon^2\sum_v|J_v|\,|u(s_v)|^2
 \geq w\int_0^L|u|^2.\label{eq:longitudinal-mass}
\end{align}

For $k\geq0$, the space
$F_k=\operatorname{span}\{\cos(m\pi s/L):0\leq m\leq k\}$ has an image
of dimension $k+1$ under $T$. The junctions add mass but no energy, so
min--max gives
\[
 \mu_k(\Omega)
 \leq\sup_{0\ne u\in F_k}
 \frac{\int_\Omega|\nabla Tu|^2}{\int_\Omega|Tu|^2}
 \leq\frac{\pi^2k^2}{L^2}.
\]

This proves \eqref{eq:chain-eigenvalues}.

For the strict count, take the cosine modes with $0\leq m<Lt$.
There are $\lceil Lt\rceil$ modes, including when $Lt$ is an integer.
Every nonzero linear combination has Rayleigh quotient below
$\pi^2t^2$, as does its lift. This proves \eqref{eq:chain-count}:
\[
 N_\Omega^N(\pi^2t^2)\geq\lceil Lt\rceil\geq Lt.
\]

\smallskip
\noindent\emph{Step 2: a subchain.}
Take $u\in H_0^1(0,L_*)$ and lift it as in Step~1 on the subchain's
rectangles and intermediate junctions. Set the lift to zero on both
terminal junctions and on the rest of the domain.

The endpoint values of $u$ vanish, and no intermediate junction has
another attachment. Thus the traces match at every interface,
and the lift belongs to $H^1(\Omega)$.

The relations \eqref{eq:longitudinal-energy}--\eqref{eq:longitudinal-mass}
hold with $L_*$ in place of $L$. The junction sum now runs only over
the intermediate junctions.

For the eigenvalue bound, take
$F_k^*=\operatorname{span}\{\sin(m\pi s/L_*):1\leq m\leq k+1\}$ for
$k\geq0$. The lifted space has dimension $k+1$, so min--max gives
\eqref{eq:subchain-eigenvalues}:
\[
 \mu_k(\Omega)
 \leq\sup_{0\ne u\in F_k^*}
 \frac{\int_\Omega|\nabla Tu|^2}{\int_\Omega|Tu|^2}
 \leq\frac{\pi^2(k+1)^2}{L_*^2}.
\]

For the strict count, there are $\lceil L_*t\rceil-1$ positive
integers $m<L_*t$. Their sine modes give a trial space of this
dimension below the energy $\pi^2t^2$.

Separately, the constant eigenfunction on $\Omega$ gives the lower
bound $1$. Taking the maximum of these two bounds proves
\eqref{eq:subchain-count}:
\[
 N_\Omega^N(\pi^2t^2)
 \geq\max\{1,\lceil L_*t\rceil-1\}
 \geq\max\{1,L_*t-1\}.
\]
\end{proof}

\begin{proposition}\label{prop:chain-criterion}
Assume that the assembly is a chain and
\begin{equation}\label{eq:chain-condition}
 BL>\frac{\pi AM}{4}.
\end{equation}
Then all positive Neumann eigenvalues satisfy $A\mu_k(\Omega)<4\pi k$.
For an assembly which contains a subchain as in
Lemma~\ref{lem:longitudinal}, the same conclusion holds provided
\begin{equation}\label{eq:subchain-condition}
 L_*^2>\pi A,\qquad B>\frac{\pi AM}{2L_*}.
\end{equation}
\end{proposition}

\begin{proof}
\smallskip
\noindent\emph{Step 1: chains.}
Let $t_0=4L/(\pi A)$. We split the energies into two ranges:
\[
 \begin{array}{ll}
 0<t<t_0:&N_\Omega^N(\pi^2t^2)\geq Lt>\pi At^2/4,\\[3pt]
 t\geq t_0:&Bt-M\geq Bt_0-M>0.
 \end{array}
\]
The first line follows from \eqref{eq:chain-count}.
The second follows from \eqref{eq:chain-condition}.
For $t\geq t_0$, insert the positive remainder $Bt-M$ into
\eqref{eq:assembly-neumann} to obtain
\[
 N_\Omega^N(\pi^2t^2)
 \geq\frac{\pi At^2}{4}+Bt-M>\frac{\pi At^2}{4}.
\]

\smallskip
\noindent\emph{Step 2: subchains.}
The condition $L_*^2>\pi A$ gives $2/L_*<2L_*/(\pi A)$.
We treat three energy ranges.

For $0<t\leq2/L_*$, the constant eigenfunction gives
\[
 N_\Omega^{N}(\pi^2t^2)\geq1>
 \frac{\pi A}{L_*^2}\geq\frac{\pi At^2}{4}.
\]
For $2/L_*<t\leq2L_*/(\pi A)$, we have
\[
 L_*t-1>\frac{L_*t}{2}\geq\frac{\pi At^2}{4},
\]
so \eqref{eq:subchain-count} gives the strict counting bound.

For $t\geq2L_*/(\pi A)$, the bound on $B$ in
\eqref{eq:subchain-condition} gives $Bt-M>0$.
Equation~\eqref{eq:assembly-neumann} then gives the strict counting bound.

Thus both the whole-chain and subchain arguments cover all positive
energies:
\[
 N_\Omega^{N}(E)>\frac{AE}{4\pi}\qquad(E>0).
\]
Evaluation at $E=\mu_k(\Omega)$ proves the stated eigenvalue inequality.
\end{proof}

A fixed finite collection of junctions gives a width criterion whose
constant is independent of the number of strips and their individual lengths.

\begin{theorem}[Uniform chain criterion]\label{thm:uniform-chain}
Fix a finite collection of regular junctions, each having one or two
ports of width $b>0$ and an open physical boundary arc. Then there is a
constant $C>0$, depending only on the collection, such that every
admissible chain built from scaled copies $\varepsilon J$ of these
junctions and from $M\geq1$ rectangular strips satisfies both strict
P\'olya inequalities as soon as
\begin{equation}\label{eq:uniform-chain}
 L>CM\varepsilon.
\end{equation}
Here $C$ does not depend on $M$ or on the individual strip lengths;
the latter need only be positive and such that the chain can be
realized admissibly.
\end{theorem}

\begin{proof}
\smallskip
\noindent\emph{Step 1: common constants.}
Choose bounds valid for the whole collection,
\[
 |J|\leq a_0,\qquad C_J^+\leq C_+,
 \qquad C_J^-\leq C_-.
\]
Choose $\tau>0$ with $\tau\leq1/b$ and $\tau\leq\sqrt{\xi_J}/\pi$ for
every reference junction $J$, where $\xi_J$ is its first mixed
eigenvalue. The parameter in \eqref{eq:assembly-tau} is therefore
at least $\tau$ for every chain under consideration.

Replacing that parameter by $\tau$ can only decrease $d_e$ and increase
$V/\tau$. This preserves the lower bound for $D$. Put
\[
 d=\min\left\{\frac13,\frac{\pi b\tau}{4}\right\},
 \qquad p=b+2C_++\frac2\tau,
 \qquad q=b+2C_-.
\]
\smallskip
\noindent\emph{Step 2: express the bounds in terms of relative width.}
A chain has $V=M+1\leq2M$ junctions. Set $r=\varepsilon M/L$.
Equations \eqref{eq:assembly-D}, \eqref{eq:assembly-B},
and \eqref{eq:assembly-area} give
\[
 D\geq L(d-pr),\qquad B\geq L(\cN-qr),\qquad
 \frac{AM}{L^2}\leq br+2a_0r^2.
\]
\smallskip
\noindent\emph{Step 3: choose the width threshold.}
Take $r_0>0$ so small that
\[
 pr_0<d,\qquad
 qr_0+\frac\pi4\bigl(br_0+2a_0r_0^2\bigr)<\cN.
\]
Set $C=1/r_0$. Condition~\eqref{eq:uniform-chain} gives $r<r_0$,
and hence $D>0$ and $BL>\pi AM/4$.
Theorem~\ref{thm:assembly} and Proposition~\ref{prop:chain-criterion}
now give, respectively,
\[
 A\lambda_k(\Omega)>4\pi k,\qquad
 A\mu_k(\Omega)<4\pi k\qquad(k\geq1).
\]
\end{proof}

If all strips have length $\ell$, the condition is
$\ell>C\varepsilon$. It allows any number of turns from the fixed
collection. Repeated junctions with holes give chains with arbitrarily
many holes. If the total length is fixed, the permitted width decreases
as the number of pieces grows.

\section{Approximation by smooth domains}\label{sec:approximation}

We prove Theorem~\ref{thm:approximation} in two stages. First fix a
packing by squares. Place a serpentine chain in each square and join
the chains by a tree outside the squares. For sufficiently many rows,
the resulting domain satisfies both strict inequalities and the
required spectral bounds.

Then vary the packing. A diagonal choice of the row numbers gives
one sequence with all the stated properties.

\subsection{A serpentine in one square}

In this section $N$ always denotes the number of rows (an integer),
not to be confused with the constant $\cN=1-\pi/4$.

Fix $0<\alpha<1$. We need a reference junction $J_\alpha$ homeomorphic
to a disk. It has two vertical ports of width $\alpha$ on $x=0$,
centered at $(0,0)$ and $(0,1)$. Its boundary is smooth away from
the closures of the ports. Each port has a rectangular collar with horizontal
boundary segments, as in Corollary~\ref{cor:smooth-networks}.
Such a junction can be found inside
\[
 [-C_0,0]\times[-\alpha/2,1+\alpha/2]
\]
for some fixed $C_0>0$.

One way to construct this junction is to start from the left half of
the annulus centered at $(0,1/2)$, with radii $(1-\alpha)/2$ and
$(1+\alpha)/2$.

Near each port endpoint, write the circular boundary arc as $y=g(x)$, where $g(0)=y_0$
and $g'(0)=0$. Replace this arc by
\[
 y=y_0+\chi(x)\bigl(g(x)-y_0\bigr),
\]
where $\chi$ is a smooth function with values in $[0,1]$ which
vanishes near $0$ and equals $1$ away from $0$.

Make each modification in a sufficiently small neighborhood of the
endpoint. The vertical bounds are then preserved, and the two boundary
arcs remain a positive distance apart. Each modified arc is horizontal
near its port and joins the unchanged circular arc smoothly.
The corners at the port endpoints disappear when the strips are attached.

We also need a cap $J_{\rm cap}$ homeomorphic to a disk.
It has vertical extent $\alpha$ and a single port of width $\alpha$.
Its boundary is smooth away from the closure of the port. The port has
a rectangular collar. Apply the same modification to a half-disk to obtain this cap.
Finally we fix $C>C_0$ so large that the horizontal extent of both
templates is strictly less than $C$.

\Needspace{22\baselineskip}
\begin{lemma}[A serpentine in one square]\label{lem:opt-serpentine}
Let the templates be as above and let $a>0$. For every sufficiently
large integer $N$ put
\begin{equation}\label{eq:serpent-parameters}
 s=\frac aN,\qquad w=\alpha s,\qquad
 y_i=\left(i-\frac12\right)s,\qquad \ell_N=a-2Cs.
\end{equation}
Then there is a smooth simply connected chain $S_N$ with
$\overline{S_N}\subset(0,a)^2$. It consists of $N$ horizontal rectangles
of length $\ell_N$ and width $w$, centered at heights $y_1,\ldots,y_N$.
The $N-1$ scaled turn junctions join successive rectangles, and
scaled caps close the two ends. Its area and total rectangle length are
\begin{align}
 |S_N|&=\alpha sN\ell_N+
 s^2\bigl((N-1)|J_\alpha|+2|J_{\rm cap}|\bigr)
       =\alpha a^2+O(N^{-1}),\label{eq:serpent-area}\\
 L_N&=N\ell_N=aN-2Ca.\label{eq:serpent-length}
\end{align}
Already the row rectangles give the distance bound
\[
 \sup_{x\in[0,a]^2}\operatorname{dist}(x,\overline{S_N})
 \leq\left(C+\frac12\right)s.
\]
The bound also holds for every domain containing the row rectangles,
including those constructed in Proposition~\ref{prop:opt-fixed-geometry}.
\end{lemma}

\begin{proof}
\smallskip
\noindent\emph{Step 1: place and join the rectangles.}
Let $N>2C$, so that $\ell_N>0$, and put
\[
 R_i=(Cs,a-Cs)\times(y_i-w/2,y_i+w/2),
 \qquad 1\leq i\leq N.
\]
Join rows $1$ and $2$, rows $3$ and $4$, and so on, at their right ends.
Join rows $2$ and $3$, rows $4$ and $5$, and so on, at their left ends.
Use translates of $sJ_\alpha$, reflected for the right turns.

Close the two free ends with copies of $sJ_{\rm cap}$.
The incidence graph is a path through $R_1,\ldots,R_N$.

\smallskip
\noindent\emph{Step 2: separation.}
Two successive turns on the same side have vertical ranges separated
by $(1-\alpha)s>0$. The top and bottom margins are
$(1-\alpha)s/2>0$.

Each reference template has horizontal extent less than $C$.
Its scaled copies therefore stay a positive distance from $x=0$ and
$x=a$. The pieces meet only at paired ports, and the closure of their
union lies in $(0,a)^2$.

\smallskip
\noindent\emph{Step 3: smoothness and topology.}
The rectangular collars make the boundary smooth across each port.

Each junction closure is a disk with disjoint closed ports.
The incidence graph is a path. Thus the gluing argument in
Section~\ref{subsec:assemblies} shows that the assembly closure is a disk.
Its interior is the required smooth simply connected domain $S_N$.

\smallskip
\noindent\emph{Step 4: area and length.}
Since the pieces have disjoint interiors, areas and lengths
add up:
\[
 |S_N|=Nw\ell_N+s^2\bigl((N-1)|J_\alpha|+2|J_{\rm cap}|\bigr),
 \qquad L_N=N\ell_N.
\]
Substituting $s=a/N$, $w=\alpha s$ and $\ell_N=a-2Cs$ gives
\eqref{eq:serpent-area} and \eqref{eq:serpent-length}.

\smallskip
\noindent\emph{Step 5: distance.}
Let $x=(x_1,x_2)\in[0,a]^2$. Choose a row height $y_i$ nearest to
$x_2$, so that $|x_2-y_i|\leq s/2$.

Let $x_1'$ be the projection of $x_1$ onto $[Cs,a-Cs]$, so that
$|x_1-x_1'|\leq Cs$. The point $(x_1',y_i)$ lies on the closed centerline
of $R_i$, hence in $\overline{S_N}$. Therefore
\[
 \operatorname{dist}(x,S_N)=\operatorname{dist}(x,\overline{S_N})
 \leq |x_1-x_1'|+|x_2-y_i|
 \leq\left(C+\frac12\right)s.
\]
Replacing the lowest left cap by an entrance leaves all row rectangles
unchanged. The distance estimate therefore remains valid in
Proposition~\ref{prop:opt-fixed-geometry}.
\end{proof}

\subsection{Packing and connecting the squares}

\begin{lemma}[Separated equal-square packing]\label{lem:opt-square-packing}
Let $U\subset\R^2$ be a nonempty bounded connected open set. Then for
every $\delta>0$ there are equal axis-parallel closed squares
$Q_1,\ldots,Q_q\subset U$, $q\geq1$, of side $a>0$, with positive pairwise
distances, such that
\begin{equation}\label{eq:square-packing}
 \left|U\setminus\bigcup_{i=1}^qQ_i\right|<\delta,
 \qquad
 d_H\left(\bigcup_{i=1}^qQ_i,\overline U\right)<\delta.
\end{equation}
\end{lemma}

\begin{proof}
For a square grid of mesh $\eta$, let $K_\eta$ be the union of its
closed cells contained in $U$. Every point of $U$ has positive distance
to $U^c$. It therefore belongs to $K_\eta$ for all sufficiently small
$\eta$. Since $U$ is bounded, dominated convergence gives
\[
 |U\setminus K_\eta|\longrightarrow0.
\]

To control the Hausdorff distance, take a finite $\delta/4$-net for
$\overline U$ consisting of points of $U$. For small $\eta$ each point
of the net lies in one of the selected cells, and so we can arrange
that
\[
 |U\setminus K_\eta|<\frac\delta2,
 \qquad d_H(K_\eta,\overline U)<\frac\delta2.
\]

Finally, shrink every selected cell about its center by one and the
same factor, close enough to one. The union $K'=\bigcup_iQ_i$ of the
shrunken cells has positive gaps between the squares, and
\[
 |K_\eta\setminus K'|<\frac\delta2,
 \qquad d_H(K_\eta,K')<\frac\delta2.
\]
The estimate \eqref{eq:square-packing} now follows from the
additivity of area and the triangle inequality.
\end{proof}

\begin{lemma}[A connecting tree with terminal stems]\label{lem:protected-tree}
Let $Q_1,\ldots,Q_q$ be the squares in \eqref{eq:square-packing},
where $q\geq2$. Write $p_i=(x_i,y_i)$ for their lower-left corners.
There are a finite embedded polygonal tree $T\subset U$, a number
$\ell>0$, and a compact polygonal set $K\subset T$ such that the
leaves of $T$ are exactly $p_1,\ldots,p_q$ and
\begin{gather}
 T\cap Q_i=\{p_i\},\qquad z_i=p_i-\ell(1,0),\label{eq:tree-leaves}\\
 T=K\cup\bigcup_{i=1}^q[z_i,p_i],\qquad
 K\cap[z_i,p_i]=\{z_i\},\label{eq:tree-stems}\\
 \dist\left(K,\bigcup_iQ_i\right)>0,\qquad
 \dist(T,\R^2\setminus U)>0.\label{eq:tree-clearance}
\end{gather}
The terminal segments $[z_i,p_i]$ are pairwise disjoint; in words,
the tree reaches each corner horizontally from the left, and apart
from these terminal segments it stays a positive distance away from
all the squares.
\end{lemma}

\begin{proof}
\smallskip
\noindent\emph{Step 1: paths outside the squares.}
The complement in $U$ of finitely many separated closed squares
compactly contained in $U$ is polygonally connected.
To see this, join two points of the complement by a polygonal path
in $U$.

Enlarge the squares slightly. Keep them disjoint and contained in $U$,
with both path endpoints outside them. Replace each portion of the
path inside an enlarged square by a boundary arc of that square.
The resulting polygonal path avoids all the original squares.

\smallskip
\noindent\emph{Step 2: protect the terminal segments.}
Choose $\rho>0$ so small that the squares $Q_i+[-4\rho,4\rho]^2$ are
pairwise disjoint and contained in $U$, and let
\[
 R_i=Q_i+[-\rho,\rho]^2,\qquad
 u_i=p_i-2\rho(1,0).
\]
By Step~1, join $u_1$ to each $u_i$ by a polygonal path in
$U\setminus\bigcup_iR_i$. Add the segments $[u_i,p_i]$.
Their union with the paths is a finite connected set $H$,
which meets each $Q_i$ only at $p_i$.

Neither the connecting paths nor the other terminal segments enter
$B(p_i,\rho)$, and hence
\[
 H\cap B(p_i,\rho)=\{p_i-t(1,0):0\leq t<\rho\}.
\]
We take $\ell=\rho/2$ and $z_i=p_i-\ell(1,0)$.

\smallskip
\noindent\emph{Step 3: select the tree.}
Subdivide $H$ at all endpoints, crossings, endpoints of collinear
overlaps, and points $z_i$. Keep each repeated subsegment only once.
This gives a finite embedded graph.

Choose a spanning tree and repeatedly prune all leaves other than the
corners $p_i$. Each edge $[z_i,p_i]$ is the only connection from $p_i$
to the other corners, so it survives both operations.
The resulting tree has exactly the leaves $p_1,\ldots,p_q$.

\smallskip
\noindent\emph{Step 4: clearances.}
Remove $(z_i,p_i]$ for each $i$ and call the remaining set $K$.
It is a finite union of closed edges disjoint from all the squares.
Compactness gives the first inequality in \eqref{eq:tree-clearance}.
The second follows because $T$ is a compact subset of $U$.
\end{proof}

\subsection{Geometry with the data fixed}

\begin{proposition}[A connected assembly for fixed data]\label{prop:opt-fixed-geometry}
Fix a packing as in Lemma~\ref{lem:opt-square-packing}, a number
$0<\alpha<1$, and turn and cap templates as above. When $q\geq2$, also
fix a tree as in Lemma~\ref{lem:protected-tree} and compatible reference
junctions at its nonterminal vertices. Each junction is a disk with
ports of width $\alpha$, smooth boundary away from the closures of
the ports, and rectangular collars. Then for all sufficiently large $N$
there is an admissible assembly $G_N$, with scale $s=a/N$ and reference
strip width $\alpha$, such that
\[
 \overline{G_N}\subset U,\qquad
 \partial G_N\in C^\infty,\qquad
 \overline{G_N}\text{ is homeomorphic to a closed disk}.
\]
If $M_N$, $V_N$, $L_N$ and $A_N$ denote, respectively, the number of
strips, the number of junctions, the total strip length and the area
of $G_N$, then
\begin{align}
 M_N&=qN+O(1),& V_N&=qN+O(1),\label{eq:approximation-counts}\\
 L_N&=qaN+O(1),& A_N&=q\alpha a^2+O(N^{-1}).\label{eq:approximation-size}
\end{align}
Moreover,
\[
 \sup_{x\in\bigcup_iQ_i}\dist(x,\overline{G_N})=O(a/N),
 \qquad L_{*,N}=aN-2Ca,
\]
where $L_{*,N}$ is the strip length of a subchain none of whose
intermediate junctions has any further attachment. Only finitely many
types of reference junctions occur, and the implied constants may
depend on all the fixed data, but not on $N$.
\end{proposition}

\begin{figure}[!htbp]
\centering
\definecolor{polyaDomain}{gray}{0.25}
\definecolor{polyaEntrance}{gray}{0.75}
\definecolor{polyaGuide}{gray}{0.50}
\begin{minipage}[t]{0.55\textwidth}
\centering
\begin{tikzpicture}[x=.93cm,y=.93cm,font=\small,
    polya strip/.style={draw=polyaDomain,line width=4.0pt,line cap=round,line join=round},
    guide/.style={draw=polyaGuide,dashed,line width=.45pt}]
  \path[use as bounding box] (-.55,-1.05) rectangle (7.45,3.75);
  \node[anchor=west,font=\bfseries\small] at (-.45,3.62) {(a) Two serpentines and an exterior path};
  \draw[polya strip,rounded corners=3pt]
    (.35,.225)--(-.40,.225)--(-.40,-.55)--(4.10,-.55)--(4.10,.225)--(4.85,.225);
  \foreach \dx/\name in {0/1,4.5/2}{
    \begin{scope}[shift={(\dx,0)}]
      \draw[guide] (0,0) rectangle (2.7,2.7);
      \node at (1.35,3.03) {$Q_{\name}$};
      \draw[polya strip]
        (.35,.225)--(2.35,.225)
        arc[start angle=-90,end angle=90,radius=.225]
        --(.35,.675)
        arc[start angle=-90,end angle=-270,radius=.225]
        --(2.35,1.125)
        arc[start angle=-90,end angle=90,radius=.225]
        --(.35,1.575)
        arc[start angle=-90,end angle=-270,radius=.225]
        --(2.35,2.025)
        arc[start angle=-90,end angle=90,radius=.225]
        --(.35,2.475);
    \end{scope}
  }
  \node[font=\footnotesize,align=center] at (3.45,-.91)
    {Only the lowest row of each square is attached.};
\end{tikzpicture}
\end{minipage}\hfill
\begin{minipage}[t]{0.43\textwidth}
\centering
\begin{tikzpicture}[x=.93cm,y=.93cm,font=\small,
    guide/.style={draw=polyaGuide,dashed,line width=.45pt},
    dim/.style={draw=black!70,line width=.45pt,<->}]
  \path[use as bounding box] (-.35,-1.05) rectangle (5.55,3.75);
  \node[anchor=west,font=\bfseries\small] at (-.30,3.62) {(b) Entrance and rectangular collar};
  \draw[guide] (1.6,3.35)--(1.6,0)--(4.7,0);
  \node[font=\footnotesize,anchor=east,text=polyaGuide] at (1.48,3.12) {$\partial Q_i$};
  \fill[polyaDomain] (-.10,.39) rectangle (4.50,.81);
  \fill[polyaEntrance] (.30,.39) rectangle (2.90,.81);
  \draw[polyaDomain,line width=.6pt] (-.10,.39)--(4.50,.39);
  \draw[polyaDomain,line width=.6pt] (-.10,.81)--(4.50,.81);
  \draw[guide] (.30,.39)--(.30,.81);
  \draw[guide] (2.90,.39)--(2.90,.81);
  \draw[polyaDomain,line width=.3906cm,line cap=butt]
    (4.50,1.80)--(2.90,1.80)
    arc[start angle=-90,end angle=-270,radius=.60]
    --(4.50,3.00);
  \draw[black!45,line width=.35pt] (4.57,.60)--(5.10,.60);
  \draw[black!45,line width=.35pt] (4.57,1.80)--(5.10,1.80);
  \draw[dim] (4.99,.60)--(4.99,1.80)
    node[midway,right] {$s$};
  \draw[dim] (-.26,.39)--(-.26,.81);
  \draw[black!45,line width=.35pt] (-.30,.39)--(-.10,.39);
  \draw[black!45,line width=.35pt] (-.30,.81)--(-.10,.81);
  \node[font=\footnotesize,anchor=south west] at (-.13,.94) {$w=\alpha s$};
  \draw[black!55,line width=.4pt] (.23,.94)--(-.26,.81);
  \draw[dim] (4.18,.81)--(4.18,1.59);
  \node[font=\scriptsize,anchor=east,fill=white,inner sep=1pt]
    at (4.03,1.20) {$(1-\alpha)s$};
  \draw[dim] (5.20,0)--(5.20,.60);
  \node[font=\scriptsize,anchor=west] at (5.28,.30) {$s/2$};
  \node[font=\footnotesize,align=center] at (1.60,-.41)
    {entrance rectangle};
  \draw[black!55,line width=.4pt] (.70,-.23)--(.70,.39);
  \node[font=\footnotesize,align=center] at (2.60,-.91)
    {Dashed port lines mark the interfaces.};
\end{tikzpicture}
\end{minipage}
\caption{Two serpentine chains joined by an exterior path (left), and
an enlarged picture of an entrance (right). The dashed squares only
indicate the positions; the gaps between the rows belong to the
exterior. The shaded entrance is the rectangle
\eqref{eq:entrance-junction}: its two vertical sides are interfaces,
and its horizontal boundary segments continue straight across them.
The rows are a distance $s$ apart and the strips have width
$w=\alpha s$. In the actual construction the turns are smoothed as
described in the text; the picture is only schematic.}
\label{fig:serpentine-entrance}
\end{figure}
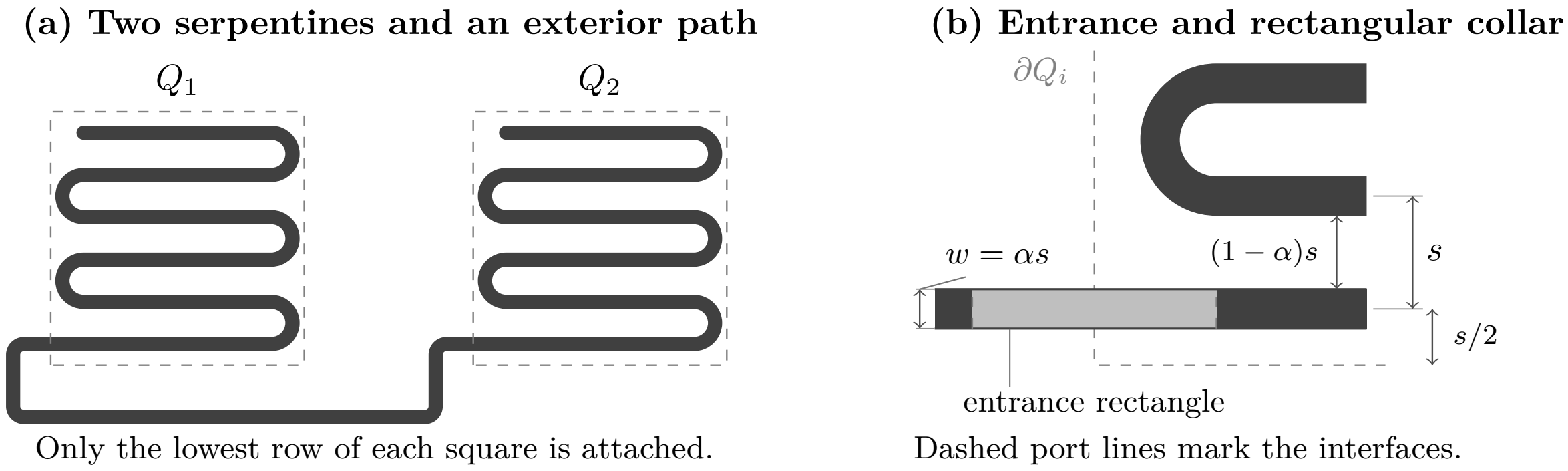

\begin{proof}
For $q=1$, use the chain of Lemma~\ref{lem:opt-serpentine}.
It has $N$ rectangles and $N+1$ junctions. The lemma gives the size
and distance bounds, smoothness, and disk topology. The whole chain
serves as the required subchain.

Assume henceforth that $q\geq2$.

\smallskip
\Needspace{10\baselineskip}
\noindent\emph{Step 1: attaching the chains.}
Place a serpentine chain in each square. Translate the tree by
$v_s=(0,s/2)$ and take $N$ large enough that $Cs<\ell$.

At the translated leaf $p_i+v_s$, replace the entrance cap by a
rectangular junction whose closure is
\begin{equation}\label{eq:entrance-junction}
 [x_i-Cs,x_i+Cs]\times
 \left[y_i+\frac{1-\alpha}{2}s,
       y_i+\frac{1+\alpha}{2}s\right].
\end{equation}
Attach its left side to the terminal stem and its right side to the
lowest row. Both sides are entire ports of width $\alpha s$.
The physical boundary continues straight across their endpoints.

At the other tree vertices, use the smoothed disk junctions of
Corollary~\ref{cor:smooth-networks}, rescaled from width $2$ to width
$\alpha$. Fix these templates before varying $N$.
Each strip along the tree has length equal to the corresponding
fixed edge length minus $O(s)$, and is positive for large $N$.

\smallskip
\noindent\emph{Step 2: clearances.}
Let
\[
 \delta_K=\dist\left(K,\bigcup_iQ_i\right)>0,
 \qquad \delta_U=\dist(T,\R^2\setminus U)>0.
\]
For the translated centerlines we have
\[
 \dist\left(K+v_s,\bigcup_iQ_i\right)\geq\delta_K-s/2,
 \qquad
 \dist(T+v_s,\R^2\setminus U)\geq\delta_U-s/2.
\]
The exterior junctions and strips lie within $O(s)$ of the translated
tree. The fixed clearances between nonincident tree parts rule out
unwanted contacts for large $N$. The protected stems separate the
entrances from the rest of the tree.

At each entrance, the next left turn lies above the entrance rectangle.
The vertical gap is
\[
 \left(y_i+\frac{3-\alpha}{2}s\right)
 -\left(y_i+\frac{1+\alpha}{2}s\right)
 =(1-\alpha)s>0.
\]
See Figure~\ref{fig:serpentine-entrance}. Since the squares and the
tree are compactly contained in $U$, we get $\overline{G_N}\subset U$.

\smallskip
\noindent\emph{Step 3: smoothness and topology.}
The rectangular collars make $\partial G_N$ smooth across the ports.

Each leaf of the exterior tree has been extended by exactly one
serpentine path. The incidence graph therefore remains a tree.
Each junction closure is a disk with disjoint closed ports, and Step~2
rules out all other contacts. The gluing argument in
Section~\ref{subsec:assemblies} gives that $\overline{G_N}$ is a disk.

Thus the gaps between rows remain connected to the exterior.
The squares only serve to position the pieces.

\smallskip
\noindent\emph{Step 4: size estimates.}
There are $qN$ row strips and $q(N-1)$ turns.
The tree pieces, entrances, and caps number $O(1)$. Hence
\[
 M_N=qN+O(1),\qquad V_N=qN+O(1),\qquad
 L_N=qN\ell_N+O(1)=qaN+O(1).
\]
The sum of the reference areas is $O(N)$, so the exact formula for
the area gives
\[
 A_N=\frac{\alpha a}{N}L_N+
          \frac{a^2}{N^2}\sum_v|J_v|
     =q\alpha a^2+O(N^{-1}).
\]
The distance estimate follows from the unchanged row rectangles,
as shown in Lemma~\ref{lem:opt-serpentine}.

\smallskip
\noindent\emph{Step 5: the long subchain.}
The serpentine in any one of the squares meets the exterior tree only
at its entrance, and none of its intermediate junctions has another
attachment. It is therefore a subchain of the required kind, of strip
length $N\ell_N=aN-2Ca$.
\end{proof}

\FloatBarrier
\Needspace{17\baselineskip}
\subsection{Spectral estimates with the data fixed}

\begin{proposition}[Spectral estimates with the data fixed]\label{prop:opt-fixed-spectrum}
Keep all the data of Proposition~\ref{prop:opt-fixed-geometry} fixed.
Then there exist $\tau>0$ and $N_0$ such that, for all $N\geq N_0$ and
all $k\geq1$,
\[
 A_N\lambda_k(G_N)>4\pi k,\qquad
 A_N\mu_k(G_N)<4\pi k.
\]
Moreover, for the same domains,
\begin{align}
 \lambda_1(G_N)&\geq\frac{\pi^2\tau^2N^2}{a^2}
             \longrightarrow\infty,\label{eq:approximation-dirichlet-limit}\\
 \mu_k(G_N)&\leq\frac{\pi^2(k+1)^2}{(aN-2Ca)^2}
             \longrightarrow0
       \quad\text{for each fixed }k\geq1.\label{eq:approximation-neumann-limit}
\end{align}
Here $N_0$ is the same for all indices.
\end{proposition}

\begin{proof}
\smallskip
\noindent\emph{Step 1: the junction errors.}
We take a lower bound which is valid for all of the finitely many
junction types:
\[
 0<\tau\leq\min\left\{\frac1\alpha,
               \min_J\frac{\sqrt{\xi_J}}\pi\right\},
 \qquad d=\min\left\{\frac13,\frac{\pi\alpha\tau}{4}\right\}>0.
\]
The parameter in \eqref{eq:assembly-tau} is at least $\tau$.

Each junction type occurs $O(N)$ times, so
\[
 \sum_v|J_v|=O(N),\qquad
 \sum_vC_{J_v}^+=O(N),\qquad
 \sum_vC_{J_v}^-=O(N).
\]
Equations \eqref{eq:assembly-D} and \eqref{eq:assembly-B}, with
$\varepsilon=a/N$, give
\begin{align*}
 D_N&\geq dL_N-\frac aN
       \left(\alpha M_N+\sum_vC_{J_v}^++\frac{V_N}{\tau}\right),\\
 B_N&=\cN L_N-\frac aN
       \left(\alpha M_N+\sum_vC_{J_v}^-\right).
\end{align*}
Thus the scaled losses stay bounded, whereas the total strip length
grows linearly in $N$:
\begin{equation}\label{eq:approximation-coefficients}
 D_N\geq dqaN-O(1),\qquad B_N=\cN qaN+O(1).
\end{equation}

\smallskip
\noindent\emph{Step 2: the Dirichlet inequality.}
For $N$ large we have $D_N>0$, and Theorem~\ref{thm:assembly} gives
both the strict Dirichlet inequality for every index and
\eqref{eq:approximation-dirichlet-limit}.

\smallskip
\noindent\emph{Step 3: the Neumann inequality.}
Lemma~\ref{lem:longitudinal}, applied to the subchain in
Proposition~\ref{prop:opt-fixed-geometry}, gives
\eqref{eq:approximation-neumann-limit}.

To obtain the inequality for every index, use the size estimates
\eqref{eq:approximation-counts}--\eqref{eq:approximation-size} and
the coefficient bound \eqref{eq:approximation-coefficients}:
\[
 \frac{L_{*,N}^2}{\pi A_N}\longrightarrow\infty,
 \qquad
 \frac{\pi A_NM_N}{2L_{*,N}}=O(1),
 \qquad B_N\longrightarrow\infty.
\]
Hence, for all large $N$,
\[
 L_{*,N}^2>\pi A_N,\qquad
 B_N>\frac{\pi A_NM_N}{2L_{*,N}}.
\]
Proposition~\ref{prop:chain-criterion} now gives
$N_{G_N}^N(E)>A_NE/(4\pi)$ for every $E>0$.
Evaluating at $E=\mu_k(G_N)$ yields
$A_N\mu_k(G_N)<4\pi k$ for every $k\geq1$.
This uses all three energy ranges in that proposition, not merely
the bound for a fixed $\mu_k$.

Choose $N_0$ to exceed the geometric threshold and the thresholds in
Steps~2 and~3. Then both P\'olya inequalities and
\eqref{eq:approximation-dirichlet-limit}--%
\eqref{eq:approximation-neumann-limit} hold on the same $G_N$ for every
$N\geq N_0$.
\end{proof}

\needspace{7\baselineskip}
\subsection{The diagonal choice}

\begin{proof}[Proof of Theorem~\ref{thm:approximation}]
\smallskip
\noindent\emph{Step 1: fix the data for each $j$.}
Fix $j\geq1$ and apply Lemma~\ref{lem:opt-square-packing} with $\delta=1/j$.
Write $K_j=\bigcup_{i=1}^{q_j}Q_{j,i}$ for the union of squares
and $a_j$ for their common side length. Choose $0<\alpha_j<1$ so that
\[
 |U\setminus K_j|<\frac1j,\qquad
 d_H(K_j,\overline U)<\frac1j,\qquad
 (1-\alpha_j)|U|<\frac1j.
\]
Fix the corresponding tree and reference templates. Only then let
$N$ vary, and denote the resulting assemblies by $G_{j,N}$.

\smallskip
\noindent\emph{Step 2: choose one row number $N_j$.}
By Propositions~\ref{prop:opt-fixed-geometry}
and~\ref{prop:opt-fixed-spectrum}, choose $N_j$ so that both strict
P\'olya inequalities hold for every index and
\begin{gather}
 \bigl||G_{j,N_j}|-q_j\alpha_ja_j^2\bigr|<\frac1j,
 \qquad
 \sup_{x\in K_j}\dist(x,\overline{G_{j,N_j}})<\frac1j,
 \label{eq:approximation-diagonal-geometry}\\
 \lambda_1(G_{j,N_j})>j,\qquad
 \mu_j(G_{j,N_j})<\frac1j.
 \label{eq:approximation-diagonal-spectrum}
\end{gather}
\smallskip
\noindent\emph{Step 3: verify the limits.}
Set $G_j=G_{j,N_j}$. These domains are smooth, simply connected,
and compactly contained in $U$. Since $q_ja_j^2=|K_j|\leq|U|$,
\begin{align*}
 |U\setminus G_j|
 &\leq (|U|-q_ja_j^2)+(1-\alpha_j)q_ja_j^2
            +\bigl|q_j\alpha_ja_j^2-|G_j|\bigr|<\frac3j,\\
 d_H(\overline{G_j},\overline U)
 &\leq d_H(K_j,\overline U)
       +\sup_{x\in K_j}\dist(x,\overline{G_j})<\frac2j.
\end{align*}
Finally, $\lambda_1(G_j)>j$, and for each fixed $k\geq1$
\[
 0\leq\mu_k(G_j)\leq\mu_j(G_j)<\frac1j\qquad(j\geq k).
\]
Thus all the assertions hold along this one sequence.
\end{proof}

The junctions and their constants may depend on $j$. They are fixed
before $N_j$ is chosen, so no bound uniform as $\alpha_j\uparrow1$
is needed.

The sequence cannot be increasing under inclusion. Otherwise Dirichlet
domain monotonicity would make $\lambda_1(G_j)$ nonincreasing,
contrary to its divergence to infinity.

\section{Rectangular frames}
\label{sec:frames}

We begin with a Neumann criterion that needs no mixed spectral gap.
Let $\Omega$ be a connected bounded Lipschitz domain, partitioned
exactly into finitely many rectangular junctions and strips.
Write $a_v,b_v$ for the side lengths of a junction, and $\ell_e,w_e$
for the length and width of a strip.

The pieces have disjoint interiors. They meet only along the designated
interfaces or at their endpoints. The longitudinal sides of every
strip lie on the physical boundary of $\Omega$. Let
\[
  A=|\Omega|,\qquad
  B=\cN\sum_e\ell_e-\sum_e w_e-\sum_v(a_v+b_v),\qquad
  M=\#\{\text{strips}\}\geq1.
\]

\begin{proposition}
\label{prop:rectangular-assembly}
If
\[
  B>\frac{M}{2}\sqrt{\pi A},
\]
then
\[
  A\mu_k(\Omega)<4\pi k\qquad(k\geq1).
\]
No lower bound on the individual strip lengths is required.
\end{proposition}

\begin{proof}
As in Lemma~\ref{lem:opt-neumann-count}, impose Dirichlet conditions
on every interface and on the boundaries of all junctions.
Keep the Neumann conditions on the longitudinal sides of the strips.
Lemma~\ref{lem:junction} gives the junction remainder constant $a_v+b_v$.

Sum the lower counts for the separated pieces.
The constant eigenfunction gives the separate lower bound $1$.
Taking the larger of these two bounds yields
\[
  N_\Omega^{N}(\pi^2t^2)
  \geq\max\left\{1,\frac{\pi At^2}{4}+Bt-M\right\}
  \qquad(t>0).
\]
For $0<t<2/\sqrt{\pi A}$, the constant bound gives $1>\pi At^2/4$.
For $t\geq2/\sqrt{\pi A}$, the hypothesis gives $Bt>M$.
The count from the separated pieces is then greater than $\pi At^2/4$.

Thus the strict Neumann inequality in~\eqref{eq:counting-target}
holds at every positive energy. Evaluation at a positive Neumann
eigenvalue proves the proposition.
\end{proof}

\begin{proof}[Proof of Theorem~\ref{thm:frames}]
After a dilation we may assume $h=1$; put $H=a+b$. The frame is then
partitioned exactly into four unit squares (the junctions) and four
strips of width one, two of length $a$ and two of length $b$, so that
\[
 M=V=4,\qquad A=2H+4,\qquad
 \sum_e\ell_e=2H,\qquad \sum_eb_e=4,\qquad
 \sum_vC_{J_v}^{\pm}=8.
\]
On a corner square, the mixed problem has Dirichlet conditions on the
two physical sides and Neumann conditions on the two sides adjacent to
the strips, and by separation of variables
\[
 \xi=\frac{\pi^2}{4}+\frac{\pi^2}{4}=\frac{\pi^2}{2},\qquad
 \tau=\frac1{\sqrt2},\qquad d_e=\frac13,
\]
\[
 D_1=\frac{2H}{3}-12-4\sqrt2\geq0
 \quad\Longleftrightarrow\quad H\geq18+6\sqrt2.
\]
By Theorem~\ref{thm:assembly}, $A\lambda_k(\Omega)>4\pi k$ for every
$k\geq1$ (the borderline case $H=18+6\sqrt2$ included).

For the Neumann eigenvalues we use
Proposition~\ref{prop:rectangular-assembly}. Here
\[
 B=2\cN H-12,\qquad
 B>2\sqrt{\pi A}
 \quad\Longleftrightarrow\quad
 \cN H-6>\sqrt{\pi(2H+4)}.
\]
With $u=H/2$, this condition becomes
\begin{equation}\label{eq:frame-numerical}
 \cN u-3>\sqrt{\pi(u+1)}.
\end{equation}
To verify this for $u\geq100$, use $\pi<22/7$ and $\cN>3/14$ and set
\[
 F(u)=\left(\frac{3u}{14}-3\right)^2-\frac{22}{7}(u+1).
\]
A direct calculation gives
\[
 F(100)=\frac{1087}{49}>0,\qquad
 F'(100)=\frac{233}{49}>0,\qquad F''(u)=\frac9{98}>0.
\]
Hence $F(u)>0$ and $3u/14-3>0$ for $u\geq100$, and therefore
\[
 \cN u-3>\frac{3u}{14}-3
   >\sqrt{\frac{22}{7}(u+1)}
   >\sqrt{\pi(u+1)}.
\]
So $H\geq200$ is sufficient, and undoing the dilation we obtain the
stated conditions on $(a+b)/h$.
\end{proof}

\begin{remark}
In the case of a square frame, the sufficient condition for the
Neumann inequalities says that the inner side length is at least
$100h$; equivalently, the ratio of the inner to the outer side length
is at least $50/51$.
\end{remark}

\section{Completion of a fixed core}
\label{sec:completion}

Let $K$ be a fixed regular connected domain. Suppose that $\partial K$
contains a straight segment $I$ of length $w>0$ lying on a supporting
line of $K$. Assume also that $\partial K\setminus I$ contains an open arc.

Attach a rectangle of width $w$ and length $L$ along $I$, on the
opposite side of the supporting line from $K$. Let $K_L$ be the interior
of the union of the two closures. The rectangle meets $K$ only along
$I$. Thus $K_L$ is a connected bounded Lipschitz domain.

Let $\xi>0$ be the first eigenvalue of $K$ with the Neumann condition
on $I$ and the Dirichlet condition on the rest of the boundary, and
let
\begin{equation}\label{eq:completion-parameters}
  \tau=\min\left\{\frac{\sqrt\xi}{\pi},\frac1w\right\},
  \qquad d=\min\left\{\frac13,\frac{\pi w\tau}{4}\right\},
\end{equation}
and choose $C_K^+,C_K^-$ as in Lemma~\ref{lem:junction}.

\Needspace{15\baselineskip}
\begin{theorem}\label{thm:completion}
Let $K$ be a regular connected domain with a straight boundary segment
$I$ of length $w>0$ lying on a supporting line. Assume that
$\partial K\setminus I$ contains an open arc, and let $K_L$ be obtained
by attaching an exterior rectangle of width $w$ and length $L$ along $I$.
With $\tau,d$ defined by \eqref{eq:completion-parameters} and
$C_K^\pm$ chosen as in Lemma~\ref{lem:junction}, both strict P\'olya
inequalities hold on $K_L$ whenever
\begin{align}
  dL&\geq w+C_K^++\frac1\tau,
  \label{eq:completion-dirichlet}\\
  \cN L-w-C_K^-&>\frac12\sqrt{\pi(|K|+wL)}.
  \label{eq:completion-neumann}
\end{align}
In particular, they hold for all sufficiently large $L$. The core is
allowed to have finitely many holes with piecewise smooth boundaries,
as long as it is regular and has a supporting boundary segment as
above.
\end{theorem}

\begin{proof}
Repeat the localization arguments of
Lemmas~\ref{lem:opt-dirichlet-localization}
and~\ref{lem:opt-neumann-count} at scale $1$.
The core is the only junction, and the attached rectangle is the only
strip. Its far end now lies on the boundary of $K_L$.

For the Dirichlet gap and upper count, replace the Dirichlet condition
at the far end by a Neumann condition. This enlarges the form domain.
For the Neumann lower count, impose a Dirichlet condition there.
This restricts the form domain.

The separated rectangle has the mixed conditions of
Lemma~\ref{lem:rectangle}.

The area and the two coefficients are now
\[
 A_L=|K|+wL,\qquad
 D_L=dL-w-C_K^+-\frac1\tau,\qquad
 B_L=\cN L-w-C_K^-.
\]
The proof of Theorem~\ref{thm:assembly}, with $M=1$, gives both
strict P\'olya inequalities under \eqref{eq:completion-dirichlet}
and~\eqref{eq:completion-neumann}.
Since $d>0$ and $\cN>0$, both coefficients grow linearly in $L$.
The square root of the area grows only like $\sqrt L$.
Both conditions therefore hold for all sufficiently large $L$.
\end{proof}

\begin{remark}
The supporting segment is essential for this construction.
The required length in Theorem~\ref{thm:completion} depends on the
core, and the added area $wL$ is not bounded uniformly.
\end{remark}

\Needspace{14\baselineskip}
\section*{Acknowledgements}
The authors thank Qing-Ming Cheng and Yong Luo for helpful discussions
and valuable suggestions.

Kefeng Liu is partially supported by the National Key Research and
Development Program of China (grant~2024YFA1013201), the Natural Science
Foundation of Chongqing (grant~CSTB2024NSCQ-LZX0040), and the Special
Project of Chongqing Municipal Science and Technology Bureau
(grant~2025CCZ015).

Liangdi Zhang is supported by the National Natural Science Foundation
of China (grant~12601080), the Scientific Research Foundation of
Chongqing University of Technology (grant~2026ZDZ012), the Youth
Project of the Science and Technology Research Program of Chongqing
Education Commission of China (grant~KJQN202601108) and the Chongqing Natural Science Foundation Innovation and Development Joint Fund (grant~CSTB2025NSCQ-LZX0059).

\section*{Declaration of AI use}
ChatGPT assisted with preliminary proof development, numerical experiments, and proofreading. Anthropic's Claude Opus was used for language editing. The authors take full responsibility for the content of this article.

\needspace{6\baselineskip}

\end{document}